\documentclass[twoside]{article}

\usepackage{hyperref} 

\usepackage{siunitx}

\usepackage{authblk} 
\usepackage{graphicx}
\usepackage{subfig}
\usepackage{booktabs}
\usepackage{amsmath}
\usepackage{amsfonts}
\usepackage{amsthm}
\usepackage{geometry}
\usepackage{amscd} 

\theoremstyle{plain}
\newtheorem{theorem}{Theorem}[section]
\theoremstyle{definition}
\newtheorem{definition}[theorem]{Definition} 
\newtheorem{example}{Example}
\theoremstyle{remark}
\newtheorem{remark}[theorem]{Remark} 

\begin{document}

\title{Explicit implementation strategy of high order edge finite elements and Schwarz preconditioning for the time-harmonic Maxwell's equations}

\date{}

\author[1]{Marcella Bonazzoli}
\author[1,2]{Victorita Dolean}
\author[3]{Fr\'ed\'eric Hecht}
\author[1]{Francesca Rapetti}
\affil[1]{Universit\'{e} C\^{o}te d'Azur, CNRS, LJAD, France. E-mail: {marcella.bonazzoli@unice.fr}, {victorita.dolean@unice.fr}, {francesca.rapetti@unice.fr}}
\affil[2]{University of Strathclyde, Glasgow, UK. E-mail: {victorita.dolean@strath.ac.uk}}
\affil[3]{UPMC Univ Paris 6, LJLL, Paris, France. E-mail: {frederic.hecht@upmc.fr}}

\renewcommand\Authfont{\normalsize}
\renewcommand\Affilfont{\itshape\small}

\maketitle

\begin{abstract}
\noindent
In this paper we focus on high order finite element approximations of the electric field combined with suitable  preconditioners, to solve the time-harmonic Maxwell's equations in waveguide configurations.
The implementation of high order curl-conforming finite elements is quite delicate, especially in the three-dimensional case. Here, we explicitly describe an implementation strategy, which has been embedded in the open source finite element software FreeFem++ (\url{http://www.freefem.org/ff++/}). In particular, we use the inverse of a generalized Vandermonde matrix to build a basis of generators in duality with the degrees of freedom, resulting in an easy-to-use but powerful interpolation operator. We carefully address the problem of applying the same Vandermonde matrix to possibly differently oriented tetrahedra of the mesh over the computational domain. 
We investigate the preconditioning for Maxwell's equations in the time-harmonic regime, which is an underdeveloped issue in the literature,  particularly for high order discretizations. In the numerical experiments, we study the effect of varying several parameters on the spectrum of the matrix preconditioned with overlapping Schwarz methods, both for 2d and 3d waveguide configurations.

\end{abstract}

\pagestyle{myheadings}
\markboth{{M.~Bonazzoli, V.~Dolean, F.~Hecht, F.~Rapetti}}{{Explicit implementation strategy of high order edge FEs and Schwarz preconditioning}}

\section{Introduction}

Developing high-speed microwave field measurement systems for wireless, medical or engineering industries is a challenging task. These systems often rely on high frequency (from $1$ to $60$\,GHz) electromagnetic wave propagation in waveguides, and the underlying mathematical model is given by Maxwell's equations.
High order finite elements methods make it possible, for a given precision, to reduce significantly the number of unknowns, and they are particularly well suited to discretize wave propagation problems since they can provide a solution with very low dispersion and dissipation errors. 
However, the resulting algebraic linear systems can be ill conditioned, so that preconditioning becomes mandatory when using iterative solvers. 

In the case of finite elements (FE) for Maxwell's equations, degrees of freedom (dofs) are associated with geometrical mesh elements other than nodes, such as edges or faces. Indeed, one needs to recognize that different physical quantities have
different properties and must be treated accordingly. Whitney finite elements are thus generally adopted \cite{Bossavit:1998:CE,Hiptmair:2002:FEC}. 
The high order version of Whitney finite elements we consider here is the one developed in \cite{Rapetti:2007:HOE,RapBos:2009:WFH} (for other possible high order finite element bases see for example \cite{AinCoy:2003:HFE, SchZag:2005:HON, Hiptmair:1999:CCF, GopGarDem:2005:NSA, ArnFalWin:2006:FEE}). 
In particular, within the family of Whitney finite elements, we address the edge elements case, which is the standard choice to describe the electric field solution of a waveguide propagation problem. 
We thus added high order edge finite elements to FreeFem++ \cite{Hecht:2012:NDF}. FreeFem++ is an open source domain specific language (DSL) specialized in solving boundary value problems (BVP) by using variational methods, and it is based on a natural transcription of the weak formulation of the considered BVP. Moreover, the user can add new finite elements to it by defining certain ingredients including an interpolation operator.
For the definition of the latter, in the high order edge elements case we need the generalized Vandermonde matrix introduced in \cite{BonRap:2015:HODu} to build a basis of generators in duality with the degrees of freedom. 
In this work, we carefully address various implementation issues, as the problem of applying the same Vandermonde matrix to possibly differently oriented simplices (triangles, tetrahedra) of the whole mesh, in order to be able to use in numerical experiments the concepts presented for just one simplex in \cite{BonRap:2015:HODu}.
Note that in FreeFem++ the basis functions are constructed locally, i.e.~in each simplex, without the need of a transformation from the reference simplex; the chosen definition of high order generators fits perfectly this local construction feature of FreeFem++ since it involves only the barycentric coordinates of the simplex. 
 
For Maxwell's equations in the \emph{time domain}, for which an implicit time discretization yields at each step a positive definite problem, there are many good solvers and preconditioners in the literature:
multigrid or auxiliary space methods, see e.g.~\cite{ReiScho:2002:AMGE,Bochev:2003:IAMG,KolVas:2008:ASAMG,BosRap:2005:PROW} for low order finite elements, \cite{Lai:2011:AMGHO} for high order ones, and Schwarz domain decomposition methods, see e.g.~\cite{Tos:2000:OSMM,DolGanGer:2009:OSM}. 
In this paper, we are interested in solving Maxwell's equations in the \emph{frequency domain}, also called the \emph{time-harmonic} Maxwell's equations: these involve the inherent difficulties of the \emph{indefinite} Helmholtz equation, which is difficult to solve for high frequencies with classical iterative methods \cite{ErnGan:2012:WDS}.
It is widely recognized that domain decomposition methods or preconditioners are key in solving efficiently Maxwell's equations in the time-harmonic regime. 

The first domain decomposition method for the time-harmonic Maxwell's equations was proposed by Despr\'es in \cite{Despres:1992:ADD}. Further improvements can be found in \cite{Collino:1997:NIM} where modified, more efficient Robin transmission conditions are used at the interfaces between subdomains. Over the last decade, optimized Schwarz methods were developed: for the first order formulation of the equations complete optimized results are known, also in the case of conductive medium \cite{DolGanGer:2009:OSM,ElBouajaji:12:OSM}, while for the second order (or curl-curl) formulation partial optimization results were obtained in various works.
Recently it has been shown that the convergence factors and the optimization process for the two formulations are the same \cite{DolGanLan:2015:ETC}. 

Nevertheless, the development of Schwarz algorithms and preconditioners for high order discretizations is still an open issue. A recent work for the non overlapping case is reported in \cite{MarsGeuz:2016:DDHO}. In the present work, we use overlapping Schwarz preconditioners based on impedance transmission conditions for high order discretizations of the curl-curl formulation of time-harmonic Maxwell's equations. 
Note that domain decomposition preconditioners are suited by construction to parallel computing, which is necessary for large scale simulations. The coupling of high order edge finite elements with domain decomposition preconditioners studied in this paper has been applied in \cite{BDRT:2017:EMF} to a large scale problem, coming from a practical application in microwave brain imaging: there, it is shown that the high order approximation of degree $2$ makes it possible to attain a given accuracy with much fewer unknowns and much less computing time than the lowest order approximation.

The paper is organized as follows. In Section~\ref{sec:formulation} we introduce the waveguide time-harmonic problem and its variational formulation. In Section~\ref{sec:discretization} we recall the definition of generators and degrees of freedom that we adopted here as high order edge FEs. Then, in Section~\ref{sec:impleFreeFem} we describe in detail the implementation issues of these FEs, the strategy developed to overcome those difficulties and the ingredients to add them as a new FE in FreeFem++. 
The overlapping Schwarz preconditioners we used are described in Section~\ref{sec:Schwarz}, followed in Section~\ref{sec:numerics} by the numerical experiments, both in two and three dimensions.

\section{The waveguide problem} 
\label{sec:formulation}

Waveguides are used to transfer electromagnetic power efficiently from one point in space, 
where an antenna is located, to another, where electronic components treat the in/out information. 
Rectangular waveguides, which are considered here,
are often used to transfer large amounts of microwave power at frequencies
greater than $2$\,GHz.  
In this section, we describe in detail the derivation of the simple but physically meaningful boundary value problem which simulates the electromagnetic wave propagation in such waveguide structures. 
To work in the frequency domain, we restrict the analysis to a time-harmonic
electromagnetic field varying with an angular frequency $\omega >0$. For all times
$t \in \mathbb{R}$, we consider the representation of the 
electric field $\boldsymbol{\mathcal{E}}$ and the magnetic field $\boldsymbol{\mathcal{H}}$ as 
\[\boldsymbol{\mathcal{E}}({\bf x},t) = \Re ({\bf E}({\bf x}) e^{{\tt i}\omega t}), 
\qquad \boldsymbol{\mathcal{H}}({\bf x},t) = \Re ({\bf H}({\bf x}) e^{{\tt i}\omega t}), \]
where ${\bf E}({\bf x})$, ${\bf H}({\bf x})$  are the complex amplitudes,
for all ${\bf x} \in \mathcal{D}$, $\mathcal{D} \subset
{\mathbb R}^3$ being the considered physical domain.
The mathematical model is thus given by the (fist order) \emph{time-harmonic Maxwell's equations}:  
\[ \nabla \times {\bf H} = {\tt i} \omega \varepsilon_{\sigma} {\bf E}, \qquad 
\nabla \times {\bf E} = -{\tt i} \omega \mu {\bf H}, 
\]
where $\mu $ is the magnetic permeability and $\varepsilon_{\sigma}$ the electric permittivity 
of the considered medium in $\mathcal{D}$.
To include dissipative effects, we work with a complex 
valued $\varepsilon_{\sigma}$, related to 
the dissipation-free electric permittivity $\varepsilon$ and the electrical conductivity $\sigma$ by
the relation 
$\varepsilon_{\sigma}= \varepsilon - {\tt i} \frac{\sigma}{\omega}$. This assumption holds
in the regions of $\mathcal{D}$ 
where the current density ${\bf J}$ is of conductive type, that is, 
${\bf J}$ and ${\bf E}$ are related by  Ohm's law
${\bf J} = \sigma {\bf E}$. Both $\varepsilon$ and $\mu$ are assumed to be 
positive, bounded functions. Expressing    
Maxwell's equations in terms of the electric field, and supposing that $\mu$ is constant,
we obtain the \emph{second order (or curl-curl) formulation}  
\begin{equation}\label{eq:MaxwellSecondOrder}
\nabla\times\left(\nabla\times \mathbf{E}\right) - {\gamma^2} \mathbf{E} 
= {\bf 0},
\end{equation}
where the (complex-valued) coefficient $\gamma$ is related to the physical parameters as follows
\[  \gamma = \sqrt{ \omega^2 \mu \varepsilon - {\tt i} \omega \mu \sigma } =
 \omega \sqrt{\mu \varepsilon_{\sigma}},  \qquad 
\varepsilon_{\sigma} = \varepsilon - {\tt i} \frac{\sigma}{\omega}.\]
Note that if $\sigma=0$, we have $\gamma = \tilde\omega$,
$\tilde\omega = \omega\sqrt{\mu \varepsilon}$ being the wavenumber.

\begin{figure} 
\centering	
\includegraphics[width=0.7\textwidth]{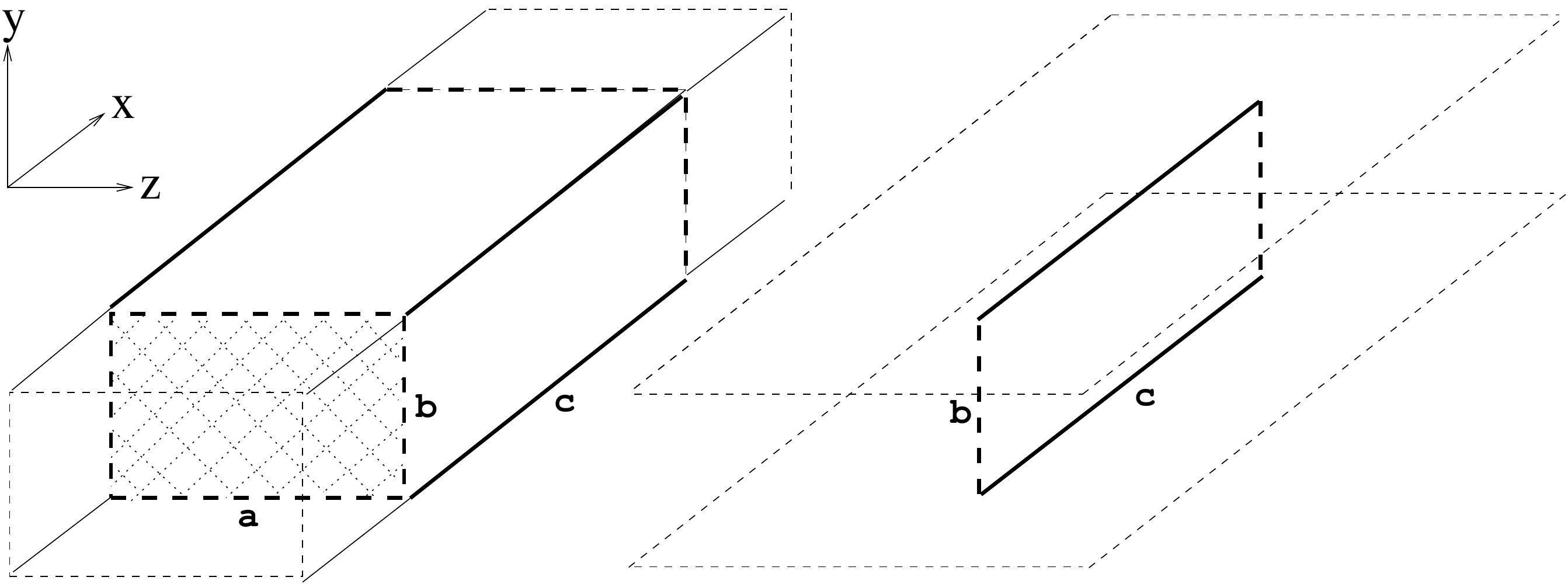}
\caption{Rectangular waveguide configuration for 3d (left) and 2d (right)  problems with 
wave propagation in the $x$-direction.
The physical domain $\mathcal{D}$ is in thin line, 
with dashed style for those boundaries that should be extended to 
infinity. The computational domain $\Omega$ is in thick line, with dashed style for those boundaries 
where suitable absorbing conditions are imposed.} 
\label{fig:dom}
\end{figure}
Equation~\eqref{eq:MaxwellSecondOrder} is to be solved in a suitable bounded section $\Omega$ of the
physical domain $\mathcal{D}$, as shown in Fig.~\ref{fig:dom}. In the 3d case, the physical domain 
$\mathcal{D} \subset \mathbb{R}^3$ is an infinite `parallelepiped'  parallel to the $x$-direction and
the computational domain is 
a  bounded section, say $\Omega=(0,{\tt X})\times (0,{\tt Y}) \times (0, {\tt Z})$ 
$= (0,{\tt c}) \times (0,{\tt b}) \times (0,{\tt a}) $  of $\mathcal{D}$.
In the 2d case, the physical domain $\mathcal{D} \subset \mathbb{R}^3$ 
is the space contained between two infinite parallel metallic plates, say   $y=0$, $y={\tt b}$, 
and all physical parameters
$\mu$, $\sigma$, $\varepsilon$ have to be assumed invariant in the $z$-direction. The computational domain
in 2d is a bounded section, say  $\Omega=(0,{\tt X})\times (0,{\tt Y})$ $=(0,{\tt c}) \times (0,{\tt b})$, 
of $\mathcal{D}$. 
In both 2d and 3d cases, the wave propagates in the $x$-direction within the domain. 

Let $\mathbf{n}$ be the unit outward normal to $\partial \Omega$.
We solve the boundary value 
problem given by equation~\eqref{eq:MaxwellSecondOrder}, with metallic boundary 
conditions
\begin{equation}\label{eq:PECbc}
\mathbf{E}\times\mathbf{n} = \mathbf{0}, \text{ on } \Gamma_{\text{w}},
\end{equation}
on the waveguide perfectly conducting walls $\Gamma_{\text{w}}= \{{\bf x} \in \partial \Omega, 
\ {\bf n}({\bf x}) \cdot {\bf e}_x = 0\}$, with ${\bf e}_x = (1,0,0)^t$,
and impedance boundary conditions
\begin{equation}\label{eq:impbc}
\begin{aligned}
& (\nabla\times\mathbf{E})\times\mathbf{n} + {\tt i} \eta \mathbf{n} \times 
(\mathbf{E} \times \mathbf{n}) =\mathbf{g}^{\text{in}}, \text{ on } \Gamma_{\text{in}} , \qquad \eta \in \mathbb{R}^+, \\
& (\nabla\times\mathbf{E})\times\mathbf{n} + {\tt i} \eta \mathbf{n} \times
(\mathbf{E} \times \mathbf{n}) = \mathbf{g}^{\text{out}}, \text{ on } \Gamma_{\text{out}},
\end{aligned}
\end{equation}
at the waveguide entrance $\Gamma_{\text{in}}=\{ {\bf x} \in \partial \Omega, 
\ {\bf n}({\bf x}) \cdot {\bf e}_x < 0 \}$, and exit $\Gamma_{\text{out}}=
\{ {\bf x} \in \partial \Omega, 
\ {\bf n}({\bf x}) \cdot {\bf e}_x > 0 \}$. 
The vectors $\mathbf{g}^{\text{in}}$, $\mathbf{g}^{\text{out}}$ depend on the incident wave.
On one hand, the impedance conditions on the artificial boundaries 
$\Gamma_{\text{in}}$, $\Gamma_{\text{out}}$ 
are absorbing boundary conditions, first order approximations of transparent boundary conditions defined to let outgoing waves pass through $\Omega$ unaffected; they
mathematically 
translate the fact that $\Omega$ is a truncated part of an infinite domain $\mathcal{D}$.
On the other hand, they simply model the fact that
the waveguide is connected to electronic components such as co-axial cables or antennas.  

To cast in the weak form the continuous problem \eqref{eq:MaxwellSecondOrder} with boundary conditions~\eqref{eq:PECbc} and \eqref{eq:impbc}, one has to multiply
\eqref{eq:MaxwellSecondOrder} by the complex conjugate of a test function $\mathbf{v}$ of a suitable
functional space $V$ and integrate by parts over the computational domain $\Omega$. More precisely,  
the weak problem reads: find $\mathbf{E} \in V$ such that
\begin{multline}
\label{eq:weakform}
\int_{\Omega} \Bigl[
(\nabla\times \mathbf{E})\cdot (\nabla\times \overline{\mathbf{v}})-\gamma^2\mathbf{E}\cdot \overline{\mathbf{v}}  \Bigr]  
+ \int_{\Gamma_{\text{in}} \cup \Gamma_{\text{out}}} {\tt i} \eta (\mathbf{E} \times \mathbf{n}) \cdot
(\overline{\mathbf{v}} \times \mathbf{n} ) \\
= \int_{\Gamma_{\text{in}}} \mathbf{g}^{\text{in}}\cdot \overline{\mathbf{v}}
+ \int_{\Gamma_{\text{out}}} \mathbf{g}^{\text{out}}\cdot \overline{\mathbf{v}} 
\quad \forall \mathbf{v}\in V,
\end{multline}
with $V=\{\mathbf{v}\in H(\text{curl},\Omega), \mathbf{v} \times \mathbf{n} = 0\text{ on }\Gamma_{\text{w}}\}$, 
where $H(\text{curl},\Omega)$ is the space of square integrable functions whose curl is also square integrable.
For a detailed discussion about existence and uniqueness of solutions we refer to \cite{Monk:2003:FEMME}.
Note that in this paper we chose the \emph{sign convention} with $e^{+{\tt i}\omega t}$ in the time-harmonic assumption, and therefore negative imaginary part in the complex valued electric permittivity $\varepsilon_{\sigma} = \varepsilon - {\tt i} \sigma/\omega$ and positive parameter $\eta$ in the impedance boundary condition~\eqref{eq:impbc}.

\section{High order edge finite elements}
\label{sec:discretization}

Consider a simplicial (triangular in 2d, tetrahedral in 3d) mesh $\mathcal{T}_h$ over $\bar{\Omega}$, where $h$ denotes the maximal diameter of simplices in $\mathcal{T}_h$.
The unknown $\mathbf{E}$ and the functional operators on it have meaningful discrete equivalents
if we work in the curl-conforming finite dimensional subspace $V_h
\subset H(\text{curl},\Omega)$  of \emph{N\'ed\'elec edge finite elements} \cite{Nedelec:1980:MFE}. 
For a simplex $T \in \mathcal{T}_h$,
the local \emph{lowest order basis functions} for the N\'ed\'elec 
curl-conforming space are associated with the oriented edges 
$ e =\{n_i,n_j\}$ of $T$ as follows
\begin{equation}\label{bf0}
{\bf w}^e = \lambda_{n_i}\nabla \lambda_{n_j} - \lambda_{n_j}\nabla\lambda_{n_i},
\end{equation}
where the $\lambda_{n_\ell}$ are the barycentric coordinates of 
a point ${\bf x} \in T$ with respect to the node $n_\ell$ of $T$ of Cartesian coordinates ${\bf x}_\ell$.
The \emph{degrees of freedom} (\emph{dofs}) over $T$ are defined as the functionals 
\[
\xi_e \colon {\bf w} \mapsto  \frac{1}{|e|} \int_e {\bf w} \cdot {\bf t}_e, \quad \forall \, e \in {\cal E}(T),
\]
where ${\bf t}_e = {\bf x}_j - {\bf x}_i$ is the tangent vector to the edge $e$, $|e| = |{\bf t}_e|$ the length of $e$ and ${\cal E}(T)$ the set of edges of $T$.  
At the lowest order, the basis functions are in \emph{duality} with the dofs, that is $\xi_e({\bf w}^{e'})=1$, resp.~$0$, if $e=e'$, resp.~if $e\ne e'$.  
As a consequence, the coefficients that define  the Galerkin projection 
${\bf E}_h$ of the field ${\bf E}$ onto $V_h$ are the \emph{circulations}  of ${\bf E}_h$
along the oriented edges $e$ of the simplicial mesh $ \mathcal{T}_h$: locally, in each $T \in \mathcal{T}_h$,
we have
\[ {\bf E}({\bf x} ) \approx {\bf E}_h({\bf x}) = 
\sum_{e \in {\cal E}(T)} c_e {\bf w}^e({\bf x}), 
\quad \forall \, {\bf x} \in T, \qquad c_e = \frac{1}{|e|} 
\int_e {\bf E}_h \cdot {\bf t}_e.\]
There are several reasons to rely on edge elements rather than on other FE
discretizations of $H(\text{curl},\Omega)$ \cite{Bossavit:1988:REE}. By construction, edge elements 
guarantee the continuity of the tangential components across inter-element interfaces,
they thus fit the continuity properties of the electric field. 
In addition, for propagation problems, edge elements are known to avoid the 
pollution of the numerical solution by spurious modes \cite{Bossavit:1990:SME, BofFerGas:1999:CMER}.

\emph{High order} curl-conforming finite elements of N\'ed\'elec type
have become  established techniques in computational electromagnetism.
Their popularity 
for wave propagation problems is due to the fact that they are characterized by  
low numerical dispersion and dissipation errors. Moreover,
at a fixed number of dofs, 
their numerical accuracy is higher. 

We adopt the high order generators of N\'ed\'elec elements 
presented in \cite{Rapetti:2007:HOE,RapBos:2009:WFH}: the definition of these generators is rather simple since it only involves the barycentric coordinates of the simplex (see also \cite{GopGarDem:2005:NSA} for previous work in this direction).
Moreover, we consider the friendly definition of the degrees of freedom and the `dualizing' Vandermonde matrix which were introduced in \cite{BonRap:2015:HODu} for all the spaces of the complex  
$
\begin{CD}
H_\text{grad} @>\nabla>>  H_\text{curl} @>\nabla \times>> H_\text{div} @>\nabla \cdot>> L^2,
\end{CD}
$.
Here, to make the presentation more accessible, we recall the definitions for the case of $H_\text{curl}$ in which we are interested, and highlight the relevant properties, also with illustrative examples. Then, in Section~\ref{sec:impleFreeFem} we describe how to deal with the delicate implementation issues of these finite elements. 

To state the definitions and further properties, we need to introduce 
multi-index notations.
A multi-index is an array ${\bf k}=(k_1, \ldots, k_{\nu})$ of $\nu$ integers
$k_i \ge 0$, and its weight $k$ is $\sum_{i=1}^{\nu} k_i$.
The set of multi-indices ${\bf k}$ with $\nu$ components and of weight
$k$ is denoted ${\cal I}(\nu,k)$. 
If $d=2,3$ is the ambient space dimension, we consider $\nu \le d+1$ and,   
given ${\bf k}\in {\cal I}(\nu,k)$, we set
$\lambda^{\bf k} = \prod_{i=1}^{\nu} \, (\lambda_{n_i})^{k_i}$,
where the $n_i$ are $\nu$ nodes of the $d+1$ nodes of $T$.
Now, in the generators definition we take $\nu=d+1$ and $k=r-1$, with $r$ the polynomial degree of the generators.

\begin{definition}[Generators]\label{wed}
{The generators for N\'ed\'elec edge element spaces $W^1_{h,r} (T)$ of \emph{degree} $r\ge1$}  
in a simplex $T \in \mathcal{T}_h$ 
are the $\lambda^{\bf k} \mathbf{w}^e$, 
with ${\bf k} \in {\cal I}(d+1,k)$, $k=r-1$ and $e \in {\cal E}(T)$. The 
$\mathbf{w}^e$ are the low order edge basis functions \eqref{bf0} (note that the polynomial degree of the $\mathbf{w}^e$ is $r=1$ and they are obtained with $k=0$).
\end{definition}

In Section 1.2 of \cite{Nedelec:1980:MFE} $W^1_{h,r}(T)$-unisolvent dofs are presented, for any $r\ge 1$ (the space $W^1_{h,r}(T)$ is indeed a discrete counterpart of $H(\text{curl},T) = \{ \mathbf{v} \in L_2(T)^3, \nabla \times  \mathbf{v} \in L_2(T)^3 \}$). 
By relying on the  generators introduced in Definition \ref{wed}, 
the functionals in \cite{Nedelec:1980:MFE} can be recast in a new more friendly form
as follows (see details in \cite{BonRap:2015:HODu}, which are inspired by \cite{Monk:2003:FEMME}).
\begin{definition}[Degrees of freedom]
\label{dofs}
For $r \ge 1$, $d=3$, the functionals
\begin{align}
& \xi_e \colon  \mathbf{w} \mapsto  \frac{1}{|e|} \int_e (\mathbf{w} \cdot \mathbf{t}_e) \, q, 
 & &\forall \, q \in \mathbb P_{r-1}(e), \  \forall \, e \in \mathcal{E}(T), \label{E2edge}  \\
& \xi_f \colon  \mathbf{w} \mapsto \frac{1}{|f|} \int_f (\mathbf{w} \cdot \mathbf{t}_{f,i})\, q,  
 & &\forall \, q \in \mathbb P_{r-2}(f), \  \forall \, f \in \mathcal{F}(T), \quad \label{E2face}\\
& & & \mathbf{t}_{f,i} \text{ two independent sides of } f, \, i = 1,2, \nonumber \\
&  \xi_T \colon \mathbf{w} \mapsto  \frac{1}{|T|}\int_T (\mathbf{w} \cdot \mathbf{t}_{T,i})\, q,   
& & \forall \, q \in \mathbb P_{r-3}(T), \label{E2vol}\\
& & & \mathbf{t}_{T,i} \text{ three independent sides of } T, \, i = 1,2,3, \nonumber
\end{align}
with $\mathcal{F}(T)$ the set of faces of $T$,
are the dofs for a function ${\bf w} \in W^1_{h,r} (T)$.
The norm of the vectors $\mathbf{t}_e, \mathbf{t}_{f,i}, \mathbf{t}_{T,i}$ is the length of the associated edge.
We say that $e, f, T$ are the \emph{supports} of the dofs $\xi_e, \xi_f, \xi_T$. 
\end{definition}

Note that for $d=2$, the dofs are given only by \eqref{E2edge} and \eqref{E2face} substituting $f$ with the triangle $T$; similarly, in the following, when $d=2$, what concerns volumes should not be taken into account and what concerns faces $f$ actually concerns the triangle $T$. 

\begin{remark}
\label{rem:polqdofs}
To make the computation of dofs easier , a \emph{convenient choice for the
polynomials} $q$ spanning the polynomial spaces over (sub)simplices $e, f,
T$ that appear in Definition~\ref{dofs} is given by suitable products
of the barycentric coordinates associated with the nodes of the
considered (sub)simplex. The space $\mathbb{P}_\rho(S)$
of polynomials of degree $\le \rho$ over a $p$-simplex $S$ (i.e.~a
simplex of dimension $1 \le p \le d$) can be generated by the products
$\lambda^{\bf k} = \prod_{i=1}^{p+1} \, (\lambda_{n_i})^{k_i}$, with
${\bf k}\in {\cal I}(p+1,\rho)$ and $n_i$ being the nodes of $S$.
\end{remark}

The \emph{classification} of dofs into edge-type, face-type, volume-type dofs can be done also for generators: volume-type generators contain (inside $\lambda^{\bf k}$ or $\mathbf{w}^e$) the barycentric coordinates w.r.t.~all the nodes of a tetrahedron $T$, face-type generators contain the ones w.r.t.~all and only the nodes of a face $f$, edge-type generators contain the ones w.r.t.~only the nodes of an edge $e$. 
Note that face-type (resp.~volume-type) generators appear for $r>1$ (resp.~$r>2$) (and the same happens for face-type and volume-type dofs).
See the explicit list of generators and dofs for the case $d=3, r=2$ in Example \ref{ex:generators}.
It turns out that dofs $\xi_e$ are $0$ on face-type and volume-type generators, and dofs $\xi_f$ are $0$ on volume-type generators.  

For the high order case ($r > 1$), the fields $\lambda^{\bf k} \mathbf{w}^e$ in Definition \ref{wed} are generators for $W^1_{h,r} (T)$, but some of the face-type or volume-type generators are \emph{linearly dependent}. The selection of generators that constitute an actual basis of $W^1_{h,r} (T)$ can be guided by the dofs in Definition~\ref{dofs}.
More precisely, as face-type (resp.~volume-type) generators keep the ones associated with the two (resp.~three) edges $e$ \emph{chosen as} the two sides $\mathbf{t}_{f,1},\mathbf{t}_{f,2}$ (resp.~three sides $\mathbf{t}_{T,1},\mathbf{t}_{T,2},\mathbf{t}_{T,3}$) of face-type dofs \eqref{E2face} (resp.~volume-type dofs \eqref{E2vol}). A convenient choice of sides is described in Subsection~\ref{subsec:implementation} and is the one adopted in Example \ref{ex:generators}. 
One can check that the total number of dofs $\xi_e, \xi_f, \xi_T$ in a simplex $T$ is equal to $\text{dim}(W^1_{h,r} (T)) = (r+d)(r+d-1)\cdots(r+2)r/(d-1)!$.

The considered basis functions are not in \emph{duality} with the dofs in Definition \ref{dofs} when $r >1$, namely,
the matrix $V$ with entries the weights $V_{ij} = \xi_i({\bf w}_{j})$, $ 1 \le i,j \le n_{\text{dofs}} = \mathrm{dim}(W^1_{h,r} (T))$ after a suitable renumbering of dofs, is not the identity matrix for $r>1$.
Duality can be re-established, if necessary, by considering new basis functions $\tilde{\mathbf{w}}_j$ built as \emph{linear combinations} of the previous basis functions with coefficients given by the entries of $V^{-1}$ \cite{BonRap:2015:HODu}. The matrix $V$ is a sort of generalized Vandermonde matrix.
Note that $V$ (and then $V^{-1}$) does not depend on the metric of the simplex $T$ for which its entries are calculated. Indeed, first of all note that dofs in Definition~\ref{dofs} are conveniently normalized. Moreover, the $\xi_i({\bf w}_{j})$ are integrals of two addends of the type $\lambda^{\bf k'} \nabla \lambda_{n_i} \cdot \mathbf{t}_e$ (here $\lambda^{\bf k'}$ gathers the products of barycentric coordinates appearing in the basis functions and in $q$, and $\mathbf{t}_e$ stands also for $\mathbf{t}_{f,i}, \mathbf{t}_{T,i}$). 
Now, we have $\nabla\lambda_{n_i} \cdot \mathbf{t}_e = -1$ if $n_i$ is the first node of $e$, $+1$ if it is its second node, $0$ if it isn't a node of $e$; so, in the end, only terms of the type $\lambda^{\bf k'}$ survive in the integral and the value of $\xi_i({\bf w}_{j})$ can be calculated using the `magic formula' (it is a classical result, see for instance \cite{Rapetti:2007:HOE}): if $S$ is a $p$-simplex, 
$\frac{1}{|S|} \int_S \prod_{i=1}^{p+1} (\lambda_{n_i})^{k_i} = p! (\prod_{i=1}^{p+1} k_i!) / (p+\sum_{i=1}^{p+1} k_i)!$. This value is clearly independent of the metric of $T$. Moreover, the entries of $V^{-1}$ turn out to be \emph{integer} numbers.
See Example \ref{ex:generators} for the case $d=3, r=2$. 

\begin{figure}
\centering
\includegraphics[width=0.25\textwidth]{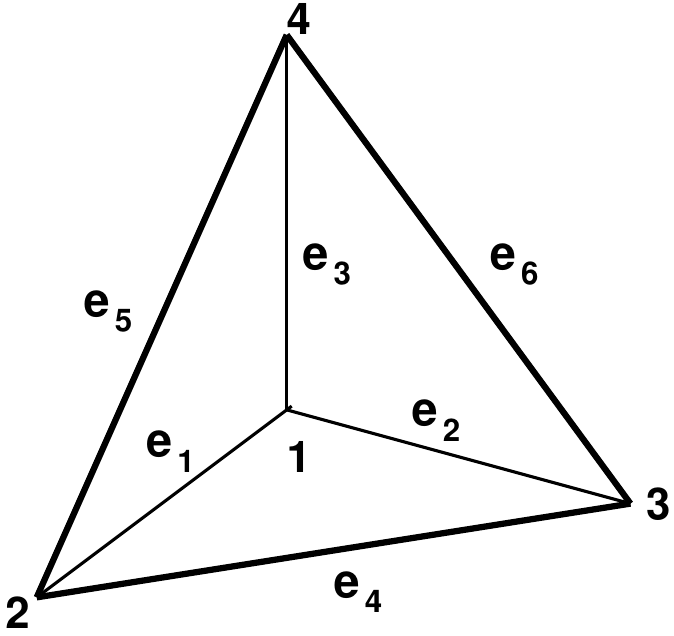}
\caption{For the tetrahedron in the figure, the edges are $e_1=\{1,2\}$, $e_2=\{1,3\}$, $e_3=\{1,4\}$, $e_4=\{2,3\}$, $e_5=\{2,4\}$, $e_6=\{3,4\}$, the faces are $f_1=\{2,3,4\}$, $f_2=\{1,3,4\}$, $f_3=\{1,2,4\}$, $f_4=\{1,2,3\}$ (note that the face $f_i$ is the one opposite the node $i$).}
\label{fig:tetNumVerEd.pdf}
\end{figure} 
\begin{example}[Generators, dofs, dualizing matrix for $d=3,$ $r=2$]
\label{ex:generators}
If the edges and the faces of a tetrahedron are numbered as in Fig.~\ref{fig:tetNumVerEd.pdf}, 
the basis functions are $\mathbf{w}_j = \lambda_{n_{r_j}} \mathbf{w}^{e_{s_j}}$, $1 \le j \le 20$,
where the $12$ edge-type basis functions have $(r_j)_{j=1}^{12} = (1,2, \,1,3, \,1,4, \,2,3, \,2,4, \,3,4)$ 
and $(s_j)_{j=1}^{12} = (1,1, \,2,2, \,3,3, \,4,4, \,5,5, \,6,6)$, 
and the $8$ face-type basis functions have $(r_j)_{j=13}^{20} = (4,3, \,4,3, \,4,2, \,3,2)$ 
and $(s_j)_{j=13}^{20} = (4,5, \,2,3, \,1,3, \,1,2)$.
Note that in order to get a basis, i.e.~a set of linearly independent generators, we have chosen to eliminate the (face-type) generators $\mathbf{w}_{21} = \lambda_{n_2} \mathbf{w}^{e_6} $, $\mathbf{w}_{22} = \lambda_{n_1} \mathbf{w}^{e_6} $,
$\mathbf{w}_{23} = \lambda_{n_1} \mathbf{w}^{e_5} $, $\mathbf{w}_{24} = \lambda_{n_1} \mathbf{w}^{e_4} $.
The corresponding edge-type dofs are:
\begin{align*}
&   \xi_1 \colon  \mathbf{w} \mapsto  \frac{1}{|e_1|} \int_{e_1} (\mathbf{w} \cdot \mathbf{t}_{e_1}) \, \lambda_{n_1} , \quad
 \xi_2 \colon  \mathbf{w} \mapsto  \frac{1}{|e_1|} \int_{e_1} (\mathbf{w} \cdot \mathbf{t}_{e_1}) \, \lambda_{n_2}, \ \dots \\
&   \xi_{11} \colon  \mathbf{w} \mapsto  \frac{1}{|e_6|} \int_{e_6} (\mathbf{w} \cdot \mathbf{t}_{e_6}) \, \lambda_{n_3} , \quad
 \xi_{12} \colon  \mathbf{w} \mapsto  \frac{1}{|e_6|} \int_{e_6} (\mathbf{w} \cdot \mathbf{t}_{e_6}) \, \lambda_{n_4} ,
\end{align*}
and the face-type dofs are:
\begin{align*}
& \xi_{13} \colon  \mathbf{w} \mapsto \frac{1}{|f_1|} \int_{f_1} (\mathbf{w} \cdot \mathbf{t}_{e_4}), \quad 
\xi_{14} \colon  \mathbf{w} \mapsto \frac{1}{|f_1|} \int_{f_1} (\mathbf{w} \cdot \mathbf{t}_{e_5}), \ \dots \\
& \xi_{19} \colon  \mathbf{w} \mapsto \frac{1}{|f_4|} \int_{f_4} (\mathbf{w} \cdot \mathbf{t}_{e_1}), \quad
\xi_{20} \colon  \mathbf{w} \mapsto \frac{1}{|f_4|} \int_{f_4} (\mathbf{w} \cdot \mathbf{t}_{e_2}). 
\end{align*}
For this ordering and choice of generators and dofs, the `dualizing' matrix $V^{-1}$ is
{
\setcounter{MaxMatrixCols}{20}
\setlength{\arraycolsep}{3pt}
\[ 
\scriptsize
V^{-1} =
\begin{bmatrix}
     4    & -2     & 0     & 0     & 0     & 0     & 0     & 0     & 0    & 0 &0  & 0     & 0     & 0     & 0     & 0     & 0     & 0     & 0     & 0  \\
    -2    & 4     & 0     & 0     & 0     & 0     & 0     & 0     & 0     & 0     & 0     & 0     & 0     & 0     & 0     & 0     & 0     & 0     & 0     & 0 \\
     0     & 0    & 4    & -2     & 0     & 0     & 0     & 0     & 0     & 0     & 0     & 0     & 0     & 0     & 0     & 0     & 0     & 0     & 0     & 0 \\
     0     & 0    & -2    & 4     & 0     & 0     & 0     & 0     & 0     & 0     & 0     & 0     & 0     & 0     & 0     & 0     & 0     & 0     & 0     & 0 \\
     0     & 0     & 0     & 0    & 4    & -2     & 0     & 0     & 0     & 0     & 0     & 0     & 0     & 0     & 0     & 0     & 0     & 0     & 0     & 0 \\
     0     & 0     & 0     & 0    & -2    & 4     & 0     & 0     & 0     & 0     & 0     & 0     & 0     & 0     & 0     & 0     & 0     & 0     & 0     & 0 \\
     0     & 0     & 0     & 0     & 0     & 0    & 4    & -2     & 0     & 0     & 0     & 0     & 0     & 0     & 0     & 0     & 0     & 0     & 0     & 0 \\
     0     & 0     & 0     & 0     & 0     & 0    & -2    & 4     & 0     & 0     & 0     & 0     & 0     & 0     & 0     & 0     & 0     & 0     & 0     & 0 \\
     0     & 0     & 0     & 0     & 0     & 0     & 0     & 0    & 4    & -2     & 0     & 0     & 0     & 0     & 0     & 0     & 0     & 0     & 0     & 0 \\
     0     & 0     & 0     & 0     & 0     & 0     & 0     & 0    & -2    & 4     & 0     & 0     & 0     & 0     & 0     & 0     & 0     & 0     & 0     & 0 \\
     0     & 0     & 0     & 0     & 0     & 0     & 0     & 0     & 0     & 0    & 4    & -2     & 0     & 0     & 0     & 0     & 0     & 0     & 0     & 0 \\
     0     & 0     & 0     & 0     & 0     & 0     & 0     & 0     & 0     & 0    & -2    & 4     & 0     & 0     & 0     & 0     & 0     & 0     & 0     & 0 \\
     0     & 0     & 0     & 0     & 0     & 0    & -4    & -2    & 2    & -2    & 2    & 4    & 8    & -4     & 0     & 0     & 0     & 0     & 0     & 0 \\
     0     & 0     & 0     & 0     & 0     & 0    & 2    & -2    & -4    & -2    & -4    & -2    & -4    & 8     & 0     & 0     & 0     & 0     & 0     & 0 \\
     0     & 0    & -4    & -2    & 2    & -2     & 0     & 0     & 0     & 0    & 2    & 4     & 0     & 0    & 8    & -4     & 0     & 0     & 0     & 0 \\
     0     & 0    & 2    & -2    & -4    & -2     & 0     & 0     & 0     & 0    & -4    & -2     & 0     & 0    & -4    & 8     & 0     & 0     & 0     & 0 \\
    -4    & -2     & 0     & 0    & 2    & -2     & 0     & 0    & 2    & 4     & 0     & 0     & 0     & 0     & 0     & 0    & 8    & -4     & 0     & 0 \\
    2    & -2     & 0     & 0    & -4    & -2     & 0     & 0    & -4    & -2     & 0     & 0     & 0     & 0     & 0     & 0    & -4    & 8     & 0     & 0 \\
    -4    & -2    & 2    & -2     & 0     & 0    & 2    & 4     & 0     & 0     & 0     & 0     & 0     & 0     & 0     & 0     & 0     & 0    & 8    & -4 \\
    2    & -2    & -4    & -2     & 0     & 0    & -4    & -2     & 0     & 0     & 0     & 0     & 0     & 0     & 0     & 0     & 0     & 0    & -4    & 8
\end{bmatrix}.
\]
}
\end{example}

\section{Implementation of high order edge finite elements in FreeFem++}
\label{sec:impleFreeFem}

In general, to \emph{add} a new finite element to FreeFem++, the user can write a C++ plugin that defines in a simplex the basis functions (and their derivatives), and an interpolation operator (which requires dofs and basis functions in \emph{duality}).
Indeed, in FreeFem++ the basis functions (and in some cases the coefficients of the interpolation operator) are constructed \emph{locally}, i.e.~in each simplex of $\mathcal{T}_h$, without the need of a transformation from the reference simplex. Note that the chosen definition of high order generators, which involves only the barycentric coordinates of the simplex, fits perfectly this local construction feature of FreeFem++. Nevertheless, the local construction should be done in such a way that the contributions coming from simplices sharing edges or faces can be then assembled properly inside the \emph{global} matrix of the FE discretization.
The strategy developed to deal with this issue for the high order edge elements is described in Subsection~\ref{subsec:implementation}.
The definition and the implementation of the interpolation operator are detailed in Subsection~\ref{sec:interpOp}. 

We added in this way the edge elements in 3d of degree $2,3$ presented before. 
The code of the C++ plugin \texttt{Element\_Mixte3d.cpp}, in which they are defined, is visible if FreeFem++ sources are downloaded (from \url{http://www.freefem.org/ff++/}) and is thus found in the folder \texttt{examples++-load}. 

\subsection{Local implementation strategy for the global assembling}
\label{subsec:implementation}
The implementation of edge finite elements is quite delicate. Indeed, basis functions and dofs are associated with the \emph{oriented} edges of mesh simplices: note that the low order $\mathbf{w}^e$ and the high order $\lambda^{\bf k} \mathbf{w}^e$ generators change sign if the orientation of the edge $e$ is reversed. 
Moreover, recall that for $r>1$, in order to get a set of linearly independent generators, we also have to \emph{choose} $2$ edges for each face $f$.
Here we wish to construct basis functions \emph{locally}, i.e.~in each simplex of $\mathcal{T}_h$, in such a way that the contributions coming from simplices sharing edges or faces could be assembled properly inside the \emph{global} matrix of the FE discretization. 
For this purpose, it is essential to orient in the \emph{same} way edges shared by simplices and to choose the \emph{same} $2$ edges for faces shared by adjacent tetrahedra. 
We have this need also to construct dofs giving the coefficients for the interpolation operator. 

\begin{figure}
\centering
\includegraphics[width=0.3\textwidth]{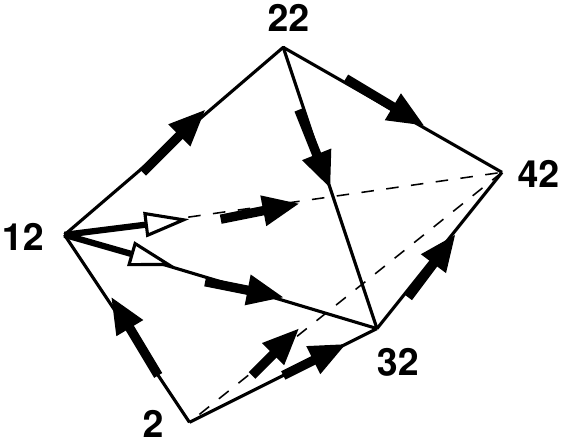}
\caption{Orientation of edges (`filled' arrows) and choice of $2$ edges (`empty' arrows) of the face shared by two adjacent tetrahedra using the numbering of mesh nodes.}
\label{fig:orientationAdjTet}
\end{figure}
This need is satisfied using the \emph{global numbers} of the mesh nodes (see Fig.~\ref{fig:orientationAdjTet}). More precisely, to orient the edges $e$ of the basis functions and the vectors $\mathbf{t}_e, \mathbf{t}_{f,i}, i=1,2, \mathbf{t}_{T,i}, i=1,2,3$ of the dofs, we go from the node with the smallest global number to the node with the biggest global number. Similarly, to choose $2$ edges per face for the face-type basis functions and dofs, we take the $2$ edges going out from the node with the smallest global number in the face (and the 1st edge goes to the node with the 2nd smallest global number, the 2nd edge goes to the node with the biggest global number in the face). 

\begin{figure}
\centering
\includegraphics[width=0.25\textwidth]{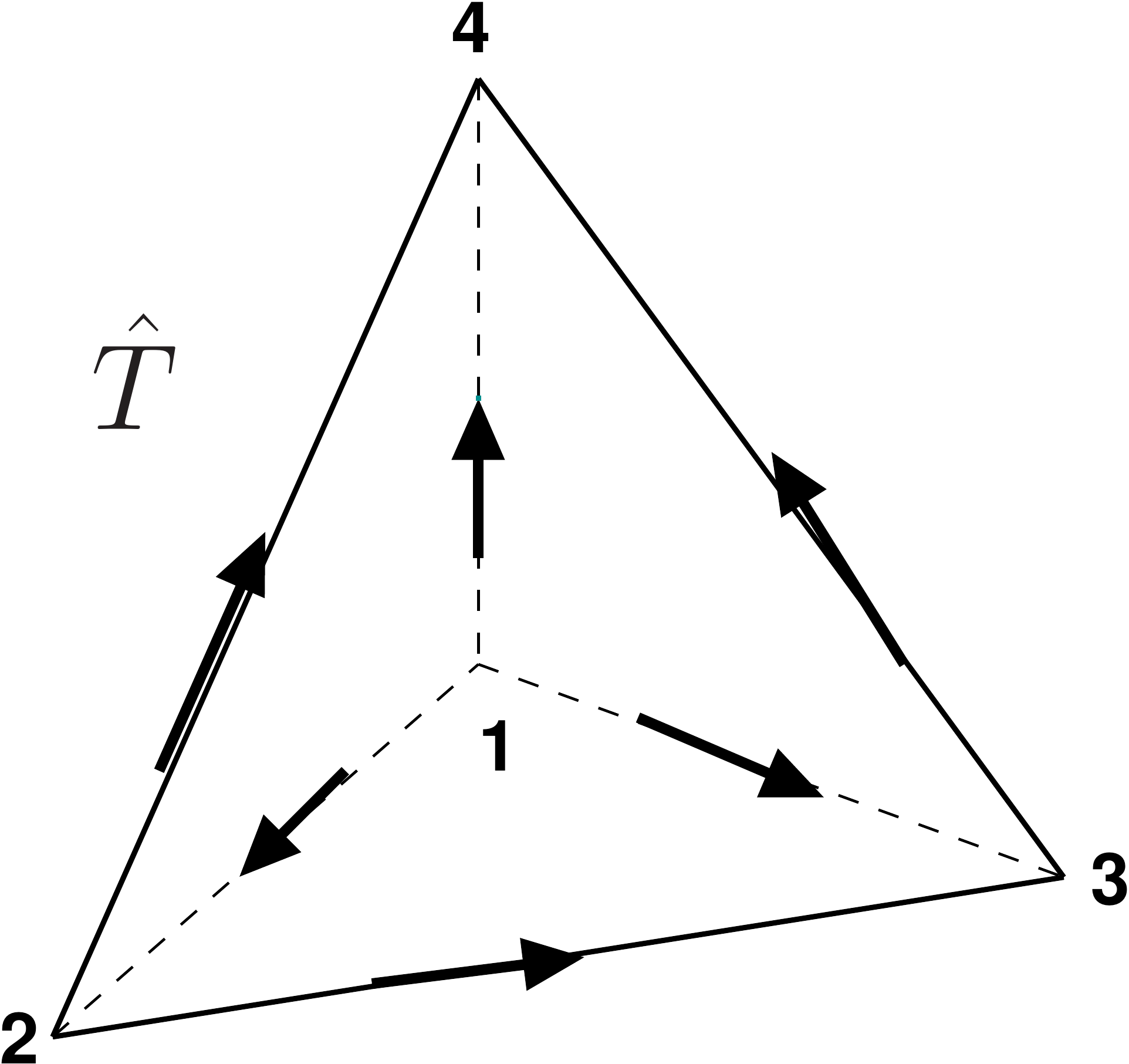}
\includegraphics[width=0.3\textwidth]{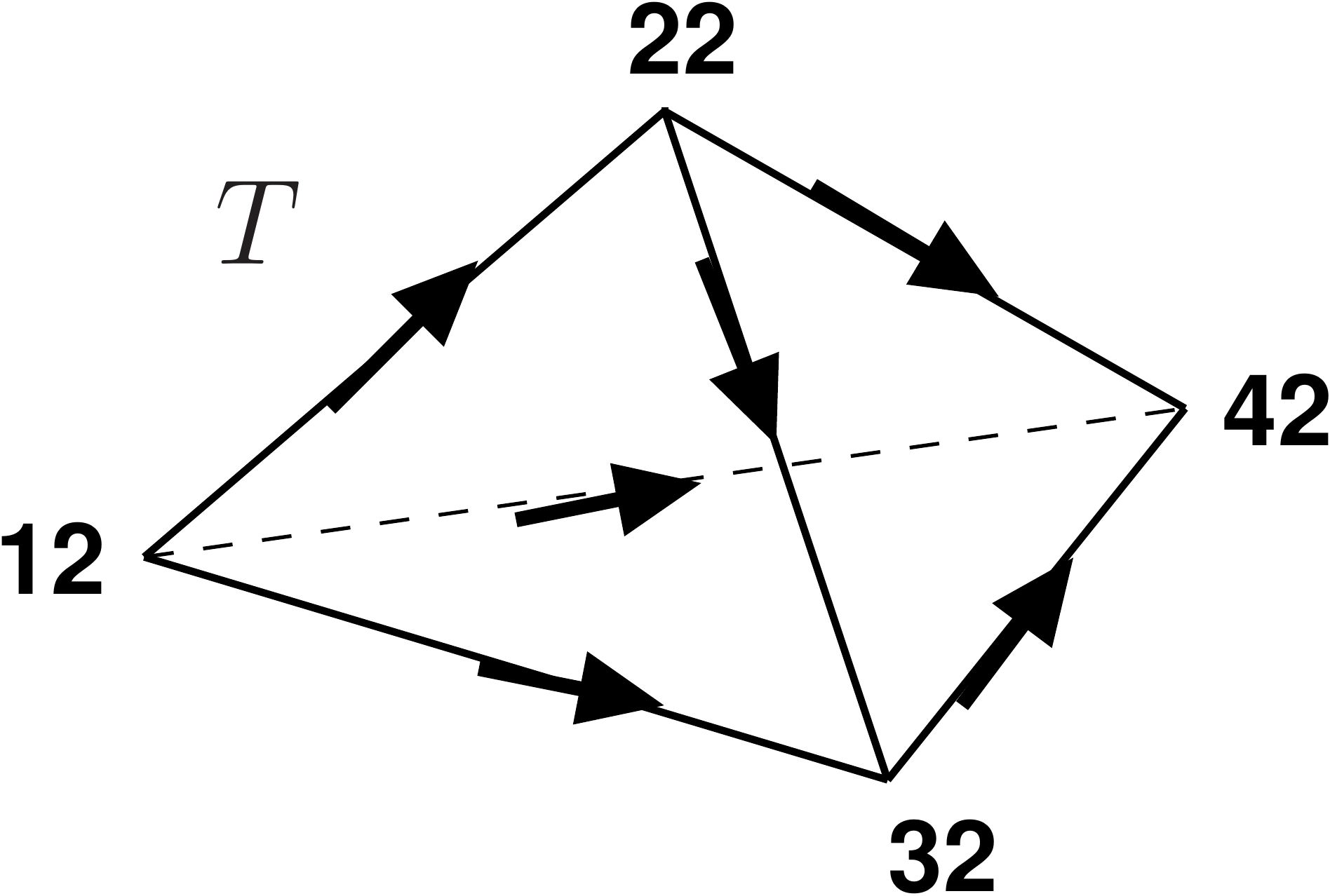}
\caption{Using global numbers to examine edges and faces, the `structure of orientation' of $T=\{12, 32, 42, 22\}$ is the one of $\hat T = \{1,2,3,4\}$ up to a rotation.}
\label{fig:orientationTThat}
\end{figure}
Moreover, when we want basis functions $\tilde{\mathbf{w}}_j$ in duality with the dofs, a \emph{second need} should be satisfied: we wish to use \emph{for all mesh simplices} $T$ the `dualizing' coefficients of the matrix $\hat V^{-1}$ calculated, once for all, for the reference simplex $\hat T$ \emph{with} a certain choice of orientation and choice of edges (recall that $V^{-1}$ already does not depend on the metric of the simplex for which it is calculated).
To be allowed to do this, it is sufficient to use the nodes \emph{global numbers} to decide the order in which the non dual $\mathbf{w}_j$ (from which we start to then get the $\tilde{\mathbf{w}}_j$) are constructed locally on $T$. 
More precisely, for the edge-type (resp.~face-type) basis functions the edges (resp.~faces) are examined in the order written in the caption of Fig.~\ref{fig:tetNumVerEd.pdf}, but replacing the nodes numbers $1,2,3,4$ with the increasing global numbers of the nodes of $T$: the 1st examined edge is from the node with the 1st smallest global number to the one with the 2nd smallest global number, the 2nd examined edge is from the node with the 1st smallest global number to the one with the 3rd smallest global number, and so on, then the 1st examined face is the one opposite the node with the smallest global number, and so on. 
Indeed, in this way the first need is respected \emph{and} the `structure of orientation' of $T$ is the one of $\hat T$ up to a rotation (see Fig.~\ref{fig:orientationTThat}): then we are allowed to use the coefficients of $\hat V^{-1}$ for the linear combinations giving the $\tilde{\mathbf{w}}_j$.

Note that in 3d (resp.~in 2d), to assemble the global linear system matrix, it is not essential which volume-type (resp.~face-type) generators are chosen since they are not shared between tetrahedra (resp.~triangles). 
On the contrary, also this choice is important when we want to use for all mesh simplices the coefficients of $\hat V^{-1}$ calculated for a simplex with a certain choice of orientation and choice of edges.

\subsubsection{Implementation of the basis functions}
To implement the strategy introduced to construct locally the basis functions $\tilde{\mathbf{w}}_j$ while respecting the two requirements just described, two \emph{permutations} can be used; note that in this paragraph the numberings start from $0$, and no more from $1$, in order to comply with the C++ plugin written for the insertion in FreeFem++ of the new FE space. 
First, to construct the non dual $\mathbf{w}_j$, we define a permutation $p_{d+1}$ of $d+1$ elements as follows: $p_{d+1}[i]$ is the local number (it takes values among $0,\dots,d$) of the node with the $i$-th smallest global number in the simplex $T$, so we can say that $p_{d+1}$ is the permutation for which the nodes of $T$ are listed with increasing global number. For instance, for the tetrahedron $T=\{12, 32, 42, 22\}$ in Fig.~\ref{fig:orientationTThat}, we have $p_4 = \{0, 3, 1, 2 \}$.
%
So, in the first step of construction of the $\mathbf{w}_j$, we replace each $\lambda_i$ appearing in their expression with $\lambda_{p_{d+1}[i]}$. In the code of the FreeFem++ plugin, the permutation $p_4$ is called \texttt{perm}.
 
Then, in the second step of construction of the $\tilde{\mathbf{w}}_j$ as linear combinations of the $\mathbf{w}_j$, we use a permutation $P_{n_{\text{dofs}}}$ of $n_{\text{dofs}} = \mathrm{dim}(W^1_{h,r} (T))$ elements to go back to the local order of edges and faces. For instance for the tetrahedron $T=\{12, 32, 42, 22\}$, the order in which edges are examined in the first step is $\{ \{12,22\},\{12,32\},\{12,42\},\{22,32\},\{22,42\},\{32,42\} \}$, while the local order of edges would be $\{ \{12,32\}, \{12,42\}, \{12,22\}, \{32,42\}, \{22,32\}, \{22,42\} \}$ (the local order is given by how the nodes of $T$ are listed); similarly, the order in which faces are examined in the first step is $\{ \{22,32,42\},$ $\{12,32,42\},$ $\{12,22,42\},$ $\{12,22,32\} \}$, while the local order of faces would be $\{ \{22,32,42\},$ $\{12,22,42\},$ $\{12,22,32\},$ $\{12,32,42\} \}$. So for this tetrahedron, if $r=2$ (for which there are $2$ basis functions for each edge and $2$ basis functions for each face, $20$ basis functions in total listed in Example \ref{ex:generators}), we have
\[ 
P_{20} = \{ 4,5,\, 0,1,\, 2,3,\, 8,9,\, 10,11,\, 6,7;\, 12,13,\, 18,19,\, 14,15,\, 16,17 \},
\] 
(note that inside each edge or face the $2$ related dofs remain ordered according to the global numbers). 
This permutation ($r=2$) is built with the following code. There, \texttt{edgesMap} corresponds to a map that associates the pair $\{ a,b \}$ of nodes of an edge $e_i$ with its number $0 \le i \le 5$; this map is rather implemented with an array defined as $\texttt{edgesMap}[(a+1)(b+1)] = i$, where $(a+1)(b+1)$ results to be unique and symmetric for a pair $(a,b)$, $0 \le a,b \le 3$, representing a tetrahedron edge.

{\footnotesize
\begin{verbatim}
int edgesMap[13] = {-1,-1,0,1,2,-1,3,-1,4,-1,-1,-1,5}; 
// static const int nvedge[6][2] = {{0,1},{0,2},{0,3},{1,2},{1,3},{2,3}};
int p20[20];
for(int i=0; i<6; ++i) // edge dofs
{
    int ii0 = Element::nvedge[i][0], ii1 = Element::nvedge[i][1];
    int i0 = perm[ii0]; int i1 = perm[ii1];
    int iEdge = edgesMap[(i0+1)*(i1+1)]; // i of the edge [i0,i1]
    p20[i*2] = iEdge*2;
    p20[i*2+1] = iEdge*2+1;
}
for(int j=0; j<4; ++j) // face dofs
{
    int jFace = perm[j];
    p20[12+j*2] = 12+jFace*2;
    p20[12+j*2+1] = 12+jFace*2+1;
}
\end{verbatim}
}
Then, we will save the linear combinations of the $\mathbf{w}_\ell$, with coefficients given by the $j$-th column of $\hat V^{-1}$ (see Example~\ref{ex:generators}), in the final basis functions $\tilde{\mathbf{w}}_{P_{20}[j]}$, thus in duality with the chosen dofs: 

{\footnotesize
\begin{verbatim}
wtilde[p20[0]] = +4*w[0]-2*w[1]-4*w[16]+2*w[17]-4*w[18]+2*w[19];
wtilde[p20[1]] = -2*w[0]+4*w[1]-2*w[16]-2*w[17]-2*w[18]-2*w[19];
wtilde[p20[2]] = +4*w[2]-2*w[3]-4*w[14]+2*w[15]+2*w[18]-4*w[19];
wtilde[p20[3]] = -2*w[2]+4*w[3]-2*w[14]-2*w[15]-2*w[18]-2*w[19];
wtilde[p20[4]] = +4*w[4]-2*w[5]+2*w[14]-4*w[15]+2*w[16]-4*w[17];
wtilde[p20[5]] = -2*w[4]+4*w[5]-2*w[14]-2*w[15]-2*w[16]-2*w[17];
wtilde[p20[6]] = +4*w[6]-2*w[7]-4*w[12]+2*w[13]+2*w[18]-4*w[19];
wtilde[p20[7]] = -2*w[6]+4*w[7]-2*w[12]-2*w[13]+4*w[18]-2*w[19];
wtilde[p20[8]] = +4*w[8]-2*w[9]+2*w[12]-4*w[13]+2*w[16]-4*w[17];
wtilde[p20[9]] = -2*w[8]+4*w[9]-2*w[12]-2*w[13]+4*w[16]-2*w[17];
wtilde[p20[10]] = +4*w[10]-2*w[11]+2*w[12]-4*w[13]+2*w[14]-4*w[15];
wtilde[p20[11]] = -2*w[10]+4*w[11]+4*w[12]-2*w[13]+4*w[14]-2*w[15];
wtilde[p20[12]] = +8*w[12]-4*w[13];
wtilde[p20[13]] = -4*w[12]+8*w[13];
wtilde[p20[14]] = +8*w[14]-4*w[15];
wtilde[p20[15]] = -4*w[14]+8*w[15];
wtilde[p20[16]] = +8*w[16]-4*w[17];
wtilde[p20[17]] = -4*w[16]+8*w[17];
wtilde[p20[18]] = +8*w[18]-4*w[19];
wtilde[p20[19]] = -4*w[18]+8*w[19];
\end{verbatim}
}

\subsection{The interpolation operator}
\label{sec:interpOp}

Duality of the basis functions with the dofs is needed in FreeFem++ to provide an \emph{interpolation operator} onto a desired FE space of a function given by its analytical expression (or of a function belonging to another FE space).
Indeed, if we define for a (vector) function $\mathbf{u}$ its finite element approximation $\mathbf{u}_h = \Pi_h(\mathbf{u})$ using the interpolation operator
\begin{equation}
\label{eq:InterpOp1}
\Pi_h \colon H(\text{curl},T) \to W^1_{h,r} (T), \quad \mathbf{u} \mapsto \mathbf{u}_h = 
\sum_{i=1}^{ n_{\text{dofs}}} c_i \tilde{\mathbf{w}}_i, \quad \text{ with } {c_i := \xi_i(\mathbf{u})},
\end{equation}
we have that, if the duality property $\xi_j(\tilde{\mathbf{w}}_i) = \delta_{ij}$ holds, then $\xi_j(\mathbf{u}_h) =  \sum_{i=1}^{ n_{\text{dofs}}} c_i \xi_j(\tilde{\mathbf{w}}_i) = c_j$.
The interpolant coefficients $c_i = \xi_i(\mathbf{u})$ are computed in FreeFem++ with suitable quadrature formulas (on edges, faces or volumes) to approximate the values of the dofs in Definition~\ref{dofs} applied to $\mathbf{u}$.

Now, denote by $g$ the whole integrand inside the dof expression, by $n_{\text{Qpts}_i}$ the number of quadrature points of the suitable quadrature formula (on a segment, triangle or tetrahedron) to compute the integral (of precision high enough  so that the integral is computed exactly when the dof is applied to a basis function), and by $\mathbf{x}_p$, $a_p$, $1 \le p \le n_{\text{QF}_i}$ the quadrature points and their weights. Then we have
\begin{equation}
\label{eq:coeffExpQuad}
c_i = \xi_i(\mathbf{u}) = \sum_{p=1}^{n_{\text{QF}_i}} a_p \, g(\mathbf{x}_{p}) = 
\sum_{p=1}^{n_{\text{QF}_i}} a_p \sum_{j=1}^d \beta_j(\mathbf{x}_{p}) \, u_j(\mathbf{x}_{p}), 
\end{equation}
where for the second equality we have factorized $g(\mathbf{x}_{p})$ in order to put in evidence the $d$ components of $\mathbf{u}$, denoted by $u_j$, $1 \le j \le d$ (see the paragraph below).

Therefore, by substituting the expression of the coefficients~\eqref{eq:coeffExpQuad} in the interpolation operator definition \eqref{eq:InterpOp1}, we have the following expression of the interpolation operator
\begin{equation}
\label{eq:InterpOp2}
\Pi_{h}(\mathbf{u}) = 
\sum_{i=1}^{n_{\text{dofs}}} \sum_{p=1}^{n_{\text{QF}_i}} \sum_{j=1}^d  a_p \, \beta_j(\mathbf{x}_{p}) \, u_j(\mathbf{x}_{p}) \, \tilde{\mathbf{w}}_i 
= \sum_{\ell=1}^{n_{\text{ind}}} \alpha_\ell \, u_{j_{\ell}}(\mathbf{x}_{p_{\ell}}) \, \tilde{\mathbf{w}}_{i_{\ell}}, 
\end{equation}
where we have set $\alpha_\ell$ equals each $a_p \, \beta_j(\mathbf{x}_{p})$ for the right triple $(i,p,j)=(i_\ell, p_\ell, j_\ell)$.
Indeed, a FreeFem++ plugin to introduce a new finite element (represented with a C++ class) should implement~\eqref{eq:InterpOp2} by specifying the quadrature points, the indices $i_\ell$ (dof indices), $p_\ell$ (quadrature point indices), $j_\ell$ (component indices), which do not depend on the simplex and are defined in the class constructor, and the coefficients $\alpha_\ell$, which can depend on the simplex (if so, which is in particular the edge elements case, the $\alpha_\ell$ are defined with the class function \texttt{set}).

\subsubsection{Interpolation operator for $d=3,$ $r=2$}

We report here the code (extracted from the plugin \texttt{Element\_Mixte3d.cpp} mentioned before) defining first the indices of~\eqref{eq:InterpOp2} for the \texttt{Edge13d} finite element, i.e.~for $d=3,$ $r=2$. There \texttt{QFe}, \texttt{QFf} are the edge, resp.~face, quadrature formulas, and \texttt{ne}=6, \texttt{nf}=4, are the number of edges, resp.~faces, of the simplex (tetrahedron); we have $n_{\text{ind}} = d \cdot \texttt{QFe.n} \cdot 2 \, \texttt{ne} + d \cdot \texttt{QFf.n} \cdot 2 \, \texttt{nf}$.
Note that in the code the numberings start from $0$, and no more from $1$. 

{\footnotesize
\begin{verbatim}
int i=0, p=0, e=0; // i is l
for(e=0; e<(Element::ne)*2; e++) // 12 edge dofs
{
  if (e%2==1) {p = p-QFe.n;} 
  // if true, the quadrature pts are the ones of the previous dof (same edge)
  for(int q=0; q<QFe.n; ++q,++p) // 2 edge quadrature pts
    for (int c=0; c<3; c++,i++) // 3 components
    {
      this->pInterpolation[i]=p; // p_l
      this->cInterpolation[i]=c; // j_l
      this->dofInterpolation[i]=e; // i_l
      this->coefInterpolation[i]=0.; // alfa_l (filled with the function set)
    }
}
for(int f=0; f<(Element::nf)*2; f++) // 8 face dofs
{
  if (f%2==1) {p = p-QFf.n;}  
  // if true, the quadrature pts are the ones of the previous dof (same face)
  for(int q=0; q<QFf.n; ++q,++p) // 3 face quadrature pts
    for (int c=0; c<3; c++,i++) // 3 components
    {
      this->pInterpolation[i]=p; // p_l
      this->cInterpolation[i]=c; // j_l
      this->dofInterpolation[i]=e+f; // i_l
      this->coefInterpolation[i]=0.; // alfa_l (filled with the function set)
    }
}
\end{verbatim}
}

Then, the coefficients $\alpha_\ell$ are defined as follows. We start by writing~\eqref{eq:coeffExpQuad} for one edge-type dof, with $e=\{n_1,n_2\}$:
\[
\begin{split}
c_i = \xi_{i}(\mathbf{u}) = \frac{1}{|e|} \int_{e} (\mathbf{u} \cdot \mathbf{t}_{e}) \, \lambda_{n_1} & =
\sum_{p=1}^{\texttt{QFe.n}} a_p \; (\mathbf{u}(\mathbf{x}_{p}) \cdot \mathbf{t}_{e}) \, \lambda_{n_1}(\mathbf{x}_{p}) \\ & = 
\sum_{p=1}^{\texttt{QFe.n}} a_p \sum_{j=1}^d u_j(\mathbf{x}_{p}) ({{x}_{n_2}}_j - {{x}_{n_1}}_j) \lambda_{n_1}(\mathbf{x}_{p})
\end{split}
\]
so $\beta_j(\mathbf{x}_{p}) = ({{x}_{n_2}}_j - {{x}_{n_1}}_j) \lambda_{n_1}(\mathbf{x}_{p})$ and
$\alpha_\ell = a_{p_\ell} \, \beta_{j_\ell}(\mathbf{x}_{p_\ell}) = ({{x}_{n_2}}_{j_\ell} - {{x}_{n_1}}_{j_\ell}) \, a_{p_\ell} \, \lambda_{n_1}(\mathbf{x}_{p_\ell})$.
Similarly for one face-type dof, with $f=\{n_1,n_2,n_3\}$, $e=\{n_1,n_2\}$:
\[
c_i = \xi_{i}(\mathbf{u}) = \frac{1}{|f|} \int_{f} (\mathbf{u} \cdot \mathbf{t}_{e}) = 
\sum_{p=1}^{\texttt{QFf.n}} a_p \; (\mathbf{u}(\mathbf{x}_{p}) \cdot \mathbf{t}_{e})  = 
\sum_{p=1}^{\texttt{QFf.n}} a_p \sum_{j=1}^d u_j(\mathbf{x}_{p}) ({{x}_{n_2}}_j - {{x}_{n_1}}_j)
\]
so $\beta_j({x}_{p}) = ({{x}_{n_2}}_j - {{x}_{n_1}}_j) $ and
$\alpha_\ell = a_{p_\ell} \, \beta_{j_\ell}(\mathbf{x}_{p_\ell}) =  ({{x}_{n_2}}_{j_\ell} - {{x}_{n_1}}_{j_\ell}) \, a_{p_\ell}$.
The code that generalizes this calculations for all the dofs is the following, extracted from the function \texttt{set} of the plugin (note that also here we have to pay particular attention to the orientation and choice issues).

{\footnotesize
\begin{verbatim}
int i=0, p=0;
for(int ee=0; ee<Element::ne; ee++) // loop on the edges
{
  R3 E=K.Edge(ee);
  int eo = K.EdgeOrientation(ee);
  if(!eo) E=-E;
    for(int edof=0; edof<2; edof++) // 2 dofs for each edge
    {
      if (edof==1) {p = p-QFe.n;}
      for(int q=0; q<QFe.n; ++q,++p)
      {
        double ll=QFe[q].x; // value of lambda_0 or lambda_1
        if( (edof+eo) == 1 ) ll = 1-ll;
        for(int c=0; c<3; c++,i++)
        {
          M.coef[i] = E[c]*QFe[q].a*ll;
        }
      }
    }
}
for(int ff=0; ff<Element::nf; ff++) // loop on the faces
{
  const Element::Vertex * fV[3] = {& K.at(Element::nvface[ff][0]), ...
  //  (one unique line with the following)
  ... & K.at(Element::nvface[ff][1]), & K.at(Element::nvface[ff][2])};
  int i0=0, i1=1, i2=2;
  if(fV[i0]>fV[i1]) Exchange(i0,i1);
    if(fV[i1]>fV[i2]) { Exchange(i1,i2);
      if(fV[i0]>fV[i1]) Exchange(i0,i1); }
  // now local numbers in the tetrahedron:
  i0 = Element::nvface[ff][i0], i1 = Element::nvface[ff][i1], ...
  ... i2 = Element::nvface[ff][i2];
  for(int fdof=0; fdof<2; ++fdof) // 2 dofs for each face
  {
    int ie0=i0, ie1 = fdof==0? i1 : i2; 
    // edge for the face dof (its endpoints local numbers)
    R3 E(K[ie0],K[ie1]);
    if (fdof==1) {p = p-QFf.n;}
    for(int q=0; q<QFf.n; ++q,++p) // loop on the 3 face quadrature pts
      for (int c=0; c<3; c++,i++) // loop on the 3 components
      {
        M.coef[i] = E[c]*QFf[q].a;
      }
  }
}
\end{verbatim}
}


\begin{example}[Using the new FEs in a FreeFem++ script]
\label{ex:FF++script}

The edge elements in 3d of degree $2,3$ can be used (since FreeFem++ version 3.44) by loading in the $\texttt{edp}$ script the plugin (\texttt{load "Element\_Mixte3d"}), and using the keywords \texttt{Edge13d}, \texttt{Edge23d} respectively. The edge elements of the lowest degree $1$ were already available and called \texttt{Edge03d}.
After generating a tetrahedral mesh \texttt{Th}, complex vector functionx \texttt{E}, \texttt{v} in, e.g., the \texttt{Edge03d} space on \texttt{Th} are declared with the commands:

{\footnotesize
\begin{verbatim}
fespace Vh(Th,Edge03d);   Vh<complex> [Ex,Ey,Ez], [vx,vy,vz];
\end{verbatim}}

\noindent Then the weak formulation~\eqref{eq:weakform} of the problem is naturally transcribed as:

{\footnotesize
\begin{verbatim}
macro Curl(ux,uy,uz) [dy(uz)-dz(uy),dz(ux)-dx(uz),dx(uy)-dy(ux)] // EOM
macro Nvec(ux,uy,uz) [uy*N.z-uz*N.y,uz*N.x-ux*N.z,ux*N.y-uy*N.x] // EOM 

problem waveguide([Ex,Ey,Ez], [vx,vy,vz], solver=sparsesolver) =
                  int3d(Th)(Curl(Ex,Ey,Ez)'*Curl(vx,vy,vz))
                - int3d(Th)(gamma^2*[Ex,Ey,Ez]'*[vx,vy,vz])
                + int2d(Th,in,out)(1i*eta*Nvec(Ex,Ey,Ez)'*Nvec(vx,vy,vz))
                - int2d(Th,in)([vx,vy,vz]'*[Gix,Giy,Giz])   
                + on(guide,Ex=0,Ey=0,Ez=0);
\end{verbatim}}

\noindent See more details in the example \texttt{waveguide.edp} available in \texttt{examples++-load} folder of every FreeFem++ distribution.
In FreeFem++ the interpolation operator is simply called with the \texttt{=} symbol: for example one can define analytical functions \texttt{func\;f1\,=\,1+x+2*y+3*z;} \texttt{func\;f2\,=\ \,-1-x-2*y+2*z;} \texttt{func\;f3\,=\ \,2-2*x+y-2*z;} and call \texttt{[Ex,Ey,Ez]=[f1,f2,f3];}.
\end{example}

\section{Overlapping Schwarz preconditioners}
\label{sec:Schwarz}

As shown numerically in \cite{Rapetti:2007:HOE}, the matrix of the linear system resulting from the described high order discretization is ill conditioned. So, when using iterative solvers (GMRES in our case), preconditioning becomes necessary, and here we choose overlapping Schwarz domain decomposition preconditioners.

Consider a decomposition of the domain $\Omega$ into $N_{\text{sub}}$ \emph{overlapping} subdomains $\Omega_s$ that consist of a union of simplices of the mesh $\mathcal{T}_h$.  
In order to describe the matrices appearing in the algebraic expression of the preconditioners, let $\mathcal{N}$ be an ordered set of the degrees of freedom of the whole domain, and let $\mathcal{N} = \bigcup_{s=1}^{N_{\text{sub}}}\mathcal{N}_s$ be its decomposition into the (non disjoint) ordered subsets corresponding to the different (overlapping) subdomains $\Omega_s$: a degree of freedom belongs to $\mathcal{N}_s$ if its support (edge, face or volume) is contained in $\Omega_s$. 
For edge finite elements, it is important to ensure that the \emph{orientation} of the degrees of freedom is the same in the domain and in the subdomains. 
Define the matrix $R_s$ as the \emph{restriction} matrix from $\Omega$ to the subdomain $\Omega_s$: it is a $\#\mathcal{N}_s \times \#\mathcal{N}$ Boolean matrix, whose $(i,j)$ entry equals $1$ if the $i$-th degree of freedom in $\mathcal{N}_s$ is the $j$-th one in $\mathcal{N}$. Note that the \emph{extension} matrix from the subdomain $\Omega_s$ to $\Omega$ is given by $R^T_s$. 
To deal with the unknowns that belong to the overlap between subdomains, define for each subdomain a $\#\mathcal{N}_s \times \#\mathcal{N}_s$ diagonal matrix $D_s$ that gives a discrete \emph{partition of unity}, i.e. 
\[
\sum_{s=1}^{N_{\text{sub}}} R^T_s D_s R_s = I.
\]
Then the \emph{Optimized Restricted Additive Schwarz} (ORAS) preconditioner can be expressed as
\begin{equation}
\label{eq:orasD}
M^{-1}_{\text{ORAS}} = \sum_{s=1}^{N_{\text{sub}}} {R}^T_s D_s A_s^{-1}R_s,
\end{equation}
where the matrices $A_s$ are the local matrices of the \emph{subproblems} with impedance boundary conditions $(\nabla\times\mathbf{E})\times\mathbf{n} + {\tt i} \tilde \omega \mathbf{n} \times 
(\mathbf{E} \times \mathbf{n})$ as transmission conditions at the interfaces between subdomains (note that in this section the term `local' refers to a subdomain and not to a mesh simplex). 
These local matrices stem from the discretization of the considered Maxwell's equation by high order finite elements introduced in the previous sections. While the term `restricted' refers to the presence of the partition of unity matrices $D_s$, the term `optimized' refers to the use of impedance boundary conditions (with parameter $\eta=\tilde\omega$ in~\eqref{eq:impbc}) as transmission conditions, proposed by Despr\'es in \cite{Despres:1992:ADD}. 
The algebraic formulation of optimized Schwarz methods, of the type of~\eqref{eq:orasD}, was introduced in \cite{StCyr:2007:OMA}.

The implementation of the partition of unity in FreeFem++ is described in \cite{2016:paperPP}: suitable piecewise linear functions $\chi_s$ giving a continuous partition of unity ($\sum_{s=1}^{N_{\text{sub}}} \chi_s=1$) are interpolated at the barycenters of the support (edge, face, volume) of each dof of the (high order) edge finite elements. This interpolation is obtained thanks to an auxiliary FreeFem++ \emph{scalar} FE space (\texttt{Edge03ds0}, \texttt{Edge13ds0}, \texttt{Edge23ds0}) that has only the interpolation operator and no basis functions, available in the plugin \texttt{Element\_Mixte3d} mentioned before.
When impedance conditions are chosen as transmission conditions at the interfaces, it is essential that not only the function $\chi_s$ but also its derivative are equal to zero on the border of the subdomain $\Omega_s$. Indeed, if this property is satisfied, the continuous version of the ORAS algorithm is equivalent to P.~L.~Lions' algorithm (see \cite{Dolean:2015:IDD}~\S2.3.2).

\section{Numerical experiments}
\label{sec:numerics}

We validate the ORAS preconditioner~\eqref{eq:orasD} for different values of physical and numerical parameters, and compare it with a symmetric variant without the partition of unity (called Optimized Additive Schwarz):
\[
M^{-1}_{\text{OAS}} = \sum_{s=1}^{N_{\text{sub}}} R^T_s A_s^{-1}R_s.
\]
The numerical experiments are performed for a waveguide configuration in $2d$ and then in $3d$.

\subsection{Results for the two-dimensional problem}

We present the results obtained for a two-dimensional waveguide with $\mathtt{c} = 0.0502$\,m, $\mathtt{b} = 0.00254$\,m, with the physical parameters: $\varepsilon = 8.85\cdot10^{-12}$\,F\,m$^{-1}$, $\mu = 1.26\cdot10^{-6}$\,H\,m$^{-1}$ and $\sigma = 0.15$\,S\,m$^{-1}$. We consider three angular frequencies $\omega_1 = 16$\,GHz, $\omega_2=32$\,GHz, and $\omega_3=64$\,GHz, varying the mesh size $h$ according to the relation $h^2\cdot\tilde{\omega}^3=2$ (in \cite{IhlBab:1995:FES} it was proved that this relation avoids pollution effects for the one-dimensional Helmholtz equation). 

Note that in 2d the function $\mathbf{E}_\text{ex} = (0, e^{-\mathtt{i}\gamma x})$ verifies the equation, the metallic boundary conditions on $\Gamma_\text{w}$, and the impedance boundary conditions on $\Gamma_{\text{in}}$, $\Gamma_{\text{out}}$ with parameter $\eta=\tilde \omega$ and $\mathbf{g}^{\text{in}} = (\mathtt{i}\gamma + \mathtt{i} \tilde\omega)\mathbf{E}_\text{ex}$ and $\mathbf{g}^{\text{out}} = (-\mathtt{i}\gamma + \mathtt{i} \tilde\omega)\mathbf{E}_\text{ex}$;
when $\sigma=0$ we get $\mathbf{g}^{\text{in}} = 2 \mathtt{i} \tilde\omega \mathbf{E}_\text{ex}$ and $\mathbf{g}^{\text{out}} = \mathbf{0}$. 
The real part of the propagation constant $-\mathtt{i}\gamma$ gives the rate at which the amplitude changes as the wave propagates, which corresponds to wave dissipation (note that if $\sigma>0$, $\Re(-\mathtt{i} \gamma)<0$, while if $\sigma=0$, $\Re(-\mathtt{i} \gamma)=0$). 
A numerical study about the order of ($h$- and $r$-) convergence with respect to the exact solution of the high order finite element method can be found in \cite{BGDR:2014:CamaProc}.

\begin{figure}
\centering
\includegraphics[width=0.65\textwidth]{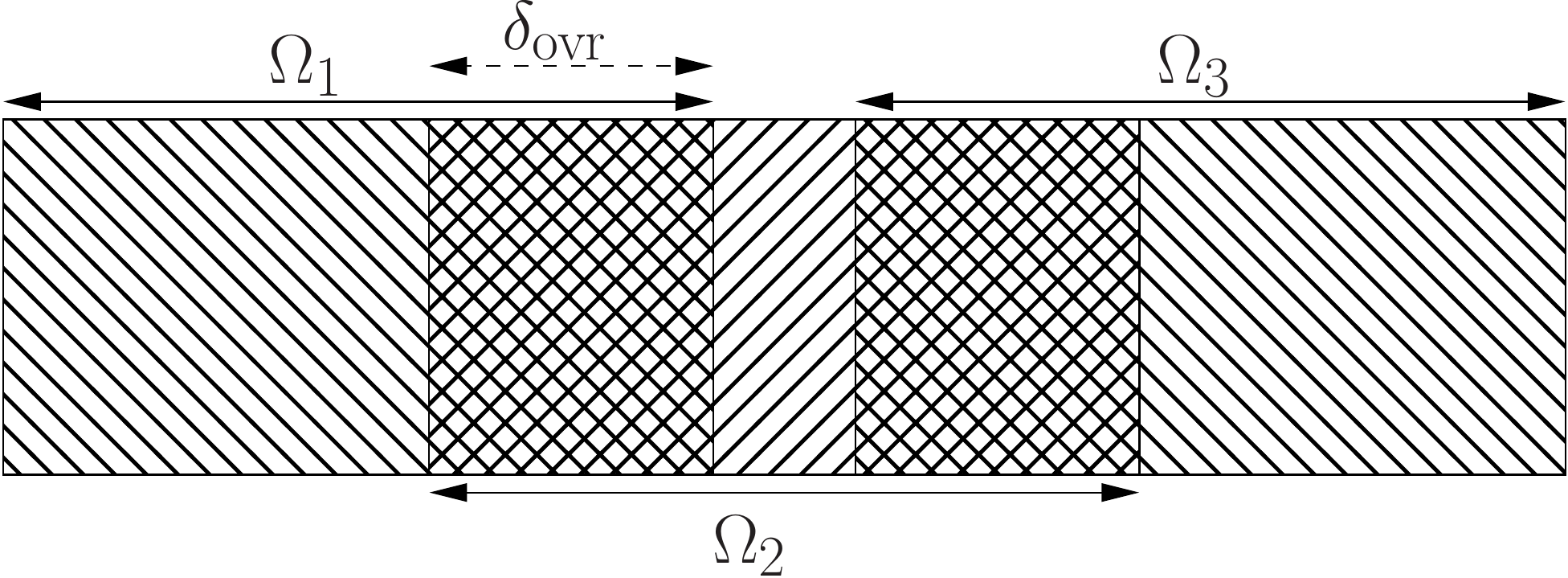}
\caption{The stripwise decomposition of the two-dimensional domain.}
\label{fig:decomp}
\end{figure}
Here we solve the linear system resulting from the finite element discretization with GMRES (with a stopping criterion based on the relative residual and a tolerance of $10^{-6}$), starting with a \emph{random} initial guess, which ensures, unlike a zero initial guess, that all frequencies are present in the error. We compare the ORAS and OAS preconditioners, taking a stripwise subdomains decomposition, along the wave propagation, as shown in Fig.~\ref{fig:decomp}. 

\begin{table}
\centering
$ \begin{array}{crrrccc}
\toprule
k & N_{\text{dofs}} & N_{\text{iterNp}} & N_{\text{iter}} & \max
\lvert \lambda-(1,0) \rvert & \#\{\lambda\in\mathbb{C}\setminus{\bar{\cal D}_1}\} & \#\{\lambda\in\partial{\cal D}_1\}  \\
\midrule
0	& 282 & 179 & 5 (10) & 1.04e{-1} (1.38e{+1}) & 0 (4) & 0 (12) \\
1	& 884 & 559 & 6 (15) & 1.05e{-1}  (1.63e{+1}) & 0 (8) & 0 (40) \\
2	& 1806 & 1138 & 6 (17) & 1.05e{-1} (1.96e{+1}) & 0 (12) & 0 (84)\\
3	& 3048 & 1946 & 6 (21) & 1.05e{-1} (8.36e{+2}) & 0 (16) & 0 (144) \\
4	& 4610 & 2950 & 6 (26) & 1.05e{-1} (1.57 e{+3}) & 0 (20) & 0 (220) \\
\bottomrule
\end{array}$
\caption{Influence of the polynomial degree $r=k+1$ on the convergence of
 ORAS(OAS) preconditioner for $\omega=\omega_2$, $N_{\text{sub}}=2$,
 $\delta_\text{ovr}=2h$.}
\label{tab:variak}
\end{table}

\begin{table}
\centering
$ \begin{array}{crrrccc}
\toprule
\omega & N_{\text{dofs}} & N_{\text{iterNp}} & N_{\text{iter}} & \max \lvert \lambda-(1,0) \rvert & \#\{\lambda\in\mathbb{C}\setminus{\bar{\cal D}_1}\} & \#\{\lambda\in\partial{\cal D}_1\} \\
\midrule
\omega_1	& 339 & 232 & 5 (11) &  2.46e{-1} (1.33 e{+1} ) & 0 (6) & 0 (45) \\
\omega_2	& 1806 & 1138 & 6 (17) & 1.05e{-1} (1.96e{+1}) & 0 (12) & 0 (84)\\
\omega_3	& 7335 & 4068 & 9 (24) & 3.03e{-1} (2.73e{+1}) & 0 (18) & 0 (123) \\
\bottomrule
\end{array}$
\caption{Influence of the angular frequency $\omega$ on the convergence of
 ORAS(OAS) preconditioner for $k=2$, $N_{\text{sub}}=2$,
 $\delta_\text{ovr}=2h$.}
\label{tab:variaw}
\end{table}

\begin{table}
\centering
$ \begin{array}{crccc}
\toprule
N_{\text{sub}} & N_{\text{iter}} & \max \lvert \lambda-(1,0) \rvert & \#\{\lambda\in\mathbb{C}\setminus{\bar{\cal D}_1}\} & \#\{\lambda\in\partial{\cal D}_1\}  \\
\midrule
2	& 6 (17) & 1.05e{-1} (1.96e{+1}) & 0 (12) & 0 (84)\\
4	& 10 (27) & 5.33e{-1} (1.96e{+1}) & 0 (38) & 0 (252) \\
8	& 19 (49) & 7.73e{-1} (1.96e{+1}) & 0 (87) & 0 (588) \\
\bottomrule
\end{array}$
\caption{Influence of the number of subdomains $N_{\text{sub}}$ on the convergence of
 ORAS(OAS) preconditioner for $k=2$, $\omega=\omega_2$,
 $\delta_\text{ovr}=2h$.}
\label{tab:variaNsub}
\end{table}

\begin{table}
\centering
$ \begin{array}{crccc}
\toprule
\delta_\text{ovr} & N_{\text{iter}} & \max \lvert \lambda-(1,0) \rvert & \#\{\lambda\in\mathbb{C}\setminus{\bar{\cal D}_1}\} & \#\{\lambda\in\partial{\cal D}_1\} \\
\midrule
1h	& 10 (20) & 1.95e{+1} (1.96e{+1}) & 3 (12) & 0 (39) \\
2h    & 6 (17) & 1.05e{-1} (1.96e{+1}) & 0 (12) & 0 (84)\\
4h	& 5 (14) & 1.06e{-1} (1.96e{+1}) & 0 (12) & 0 (174) \\			
\bottomrule
\end{array}$
\caption{Influence of the overlap size $\delta_\text{ovr}$ on the convergence of
 ORAS(OAS) preconditioner for $k=2$, $\omega=\omega_2$,
 $N_{\text{sub}}=2$.}
\label{tab:variaovr}
\end{table}

\begin{figure}
\centering
\subfloat[][$k=0$]
  {\includegraphics[width=.45\textwidth]{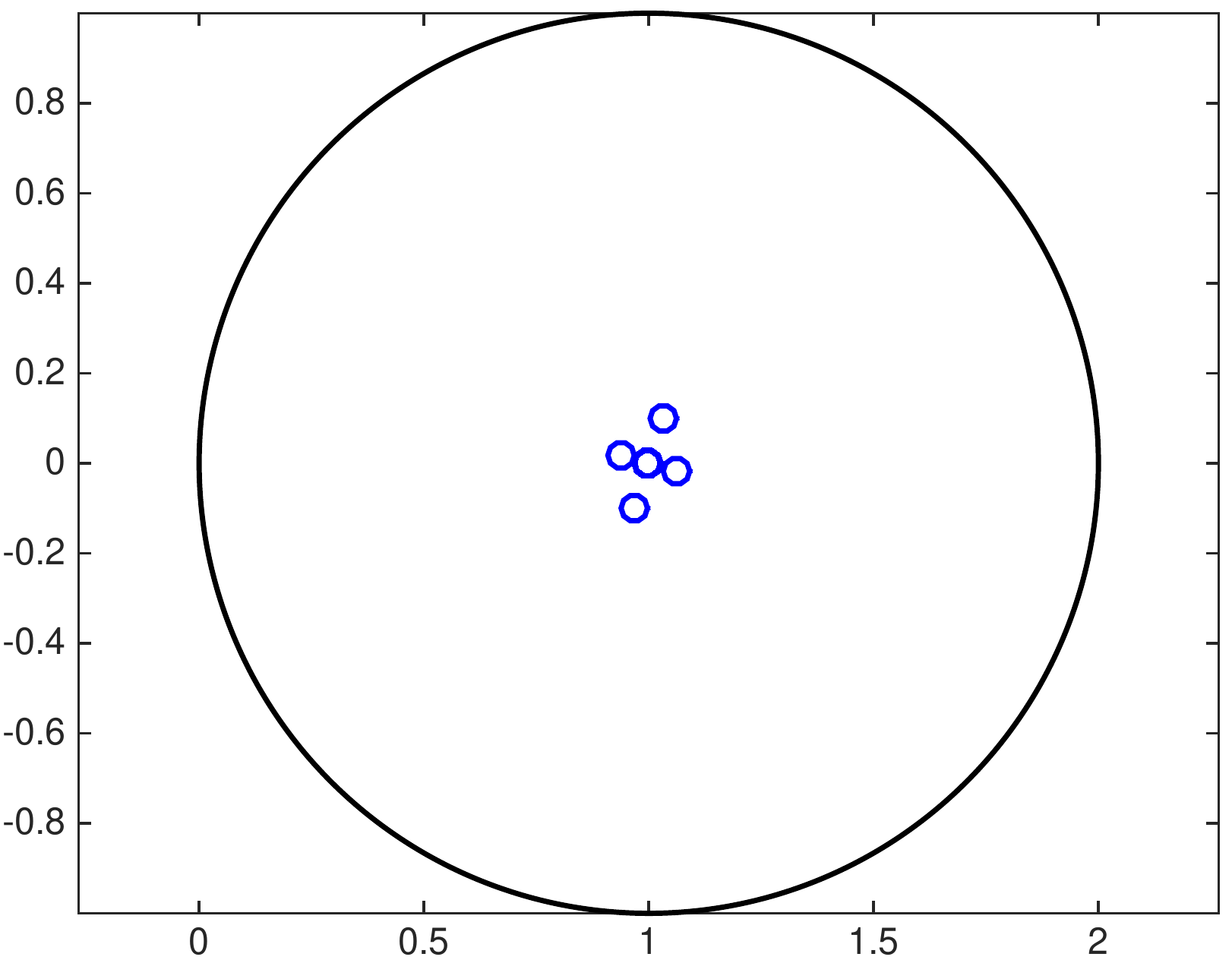}} \quad
\subfloat[][$k=1,2,3$]
  {\includegraphics[width=.45\textwidth]{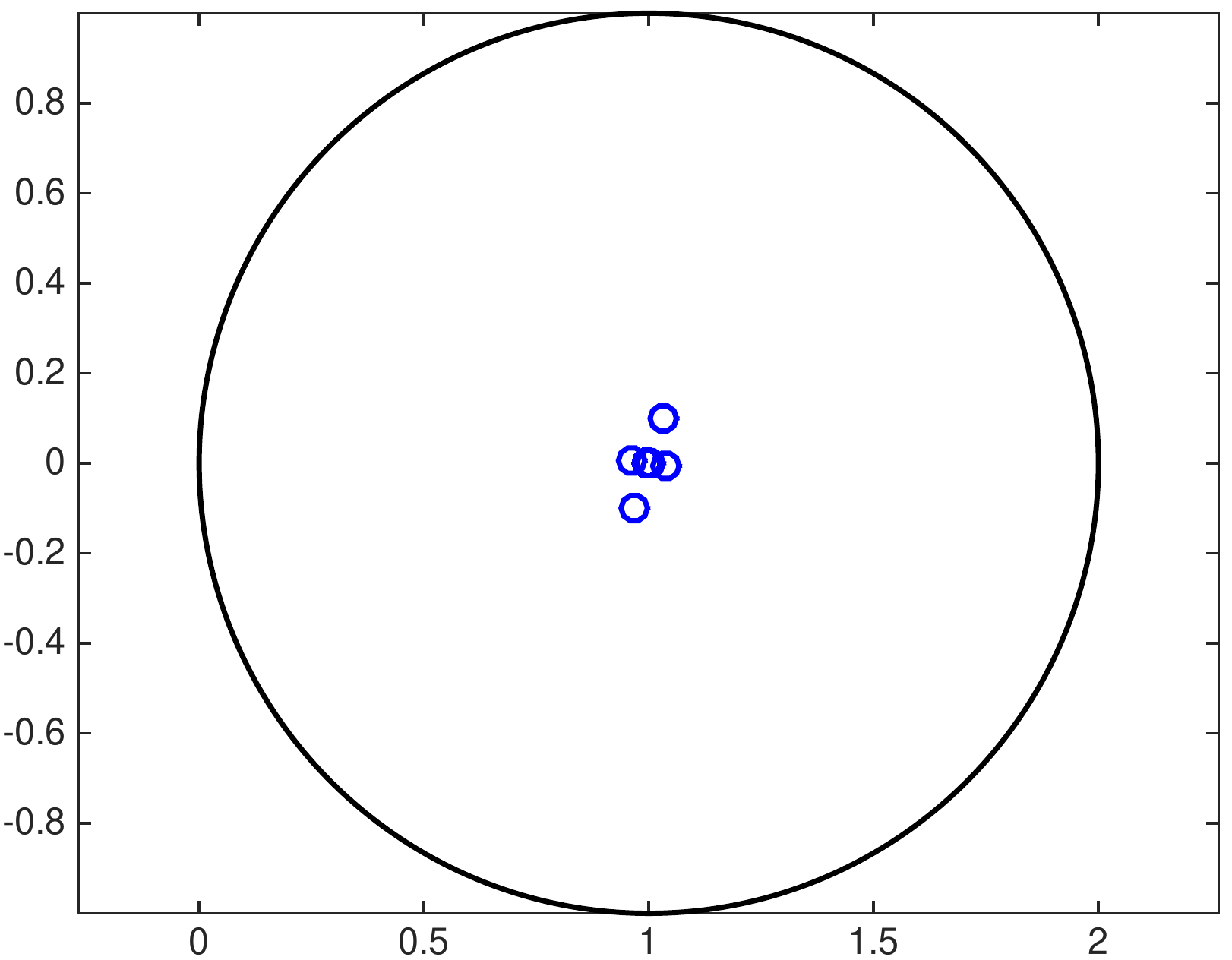}} 
\caption{Influence of the polynomial degree $r=k+1$ on the spectrum of the ORAS-preconditioned matrix for $\omega=\omega_2$, $N_{\text{sub}}=2$, $\delta_\text{ovr}=2h$.}
\label{fig:spectraRASvariak}
\end{figure}

\begin{figure}
\centering
\subfloat[][$k=0$]
  {\includegraphics[width=.45\textwidth]{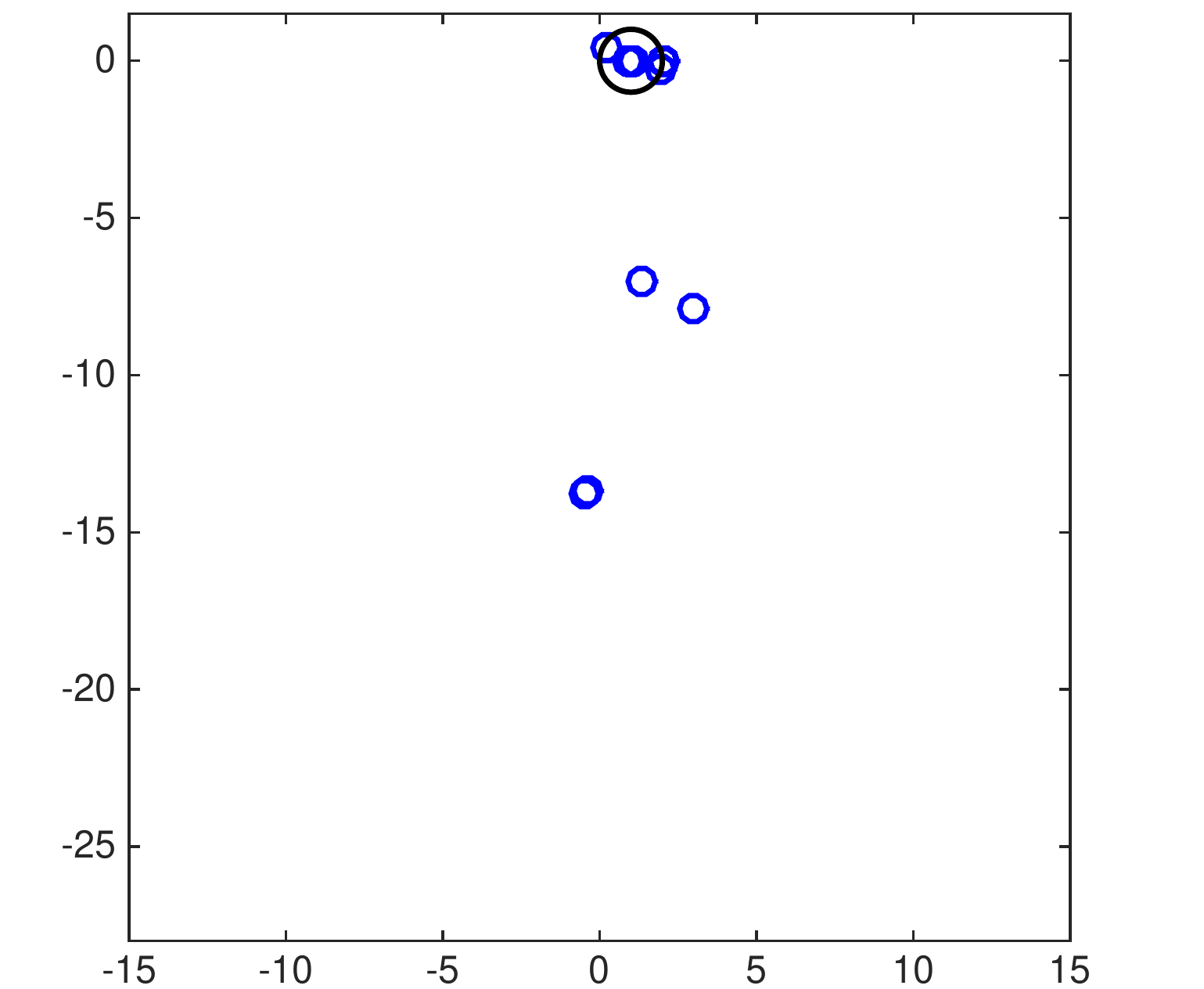}} \quad
\subfloat[][$k=1$]
  {\includegraphics[width=.45\textwidth]{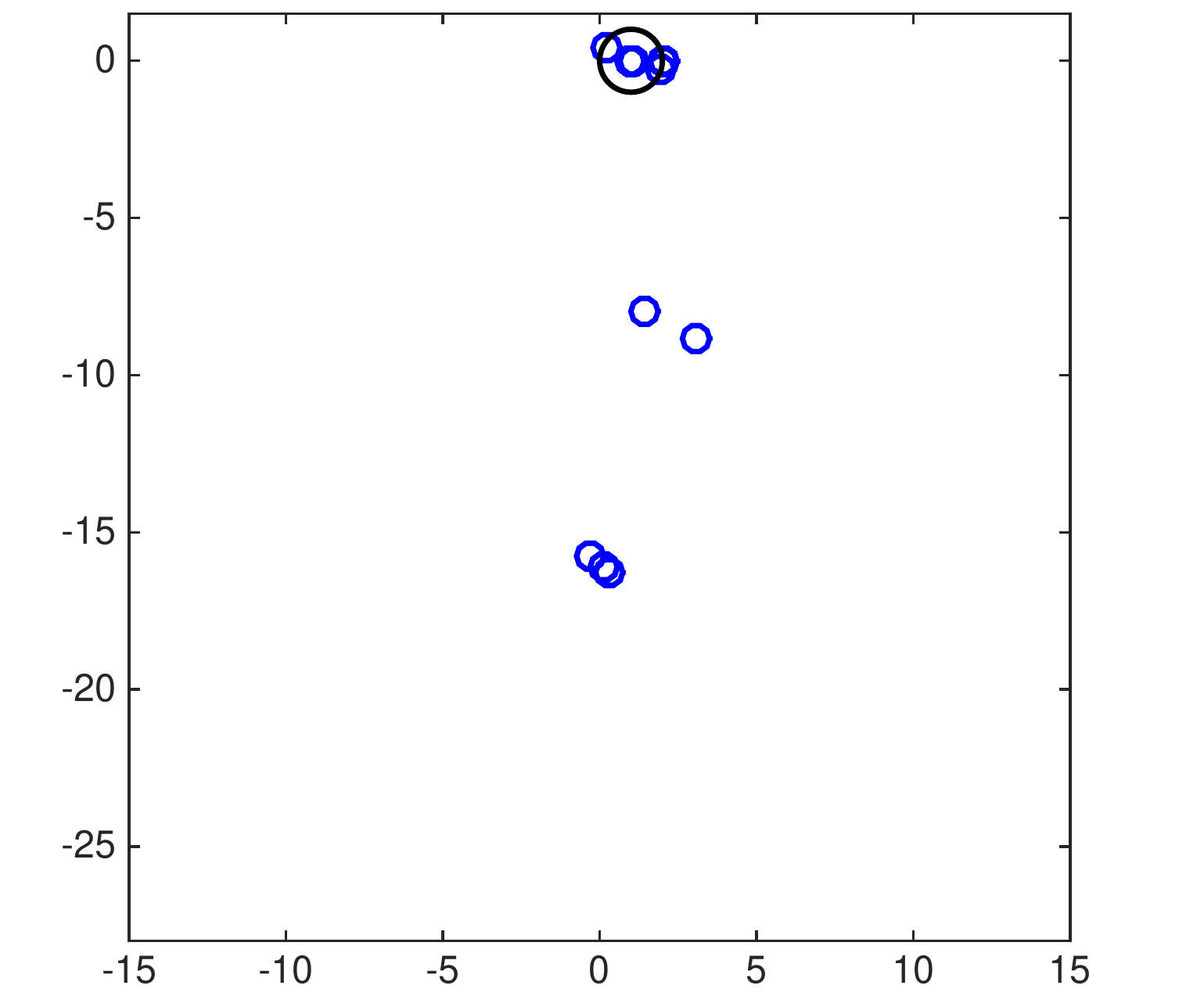}} \\
\subfloat[][$k=2$]
  {\includegraphics[width=.45\textwidth]{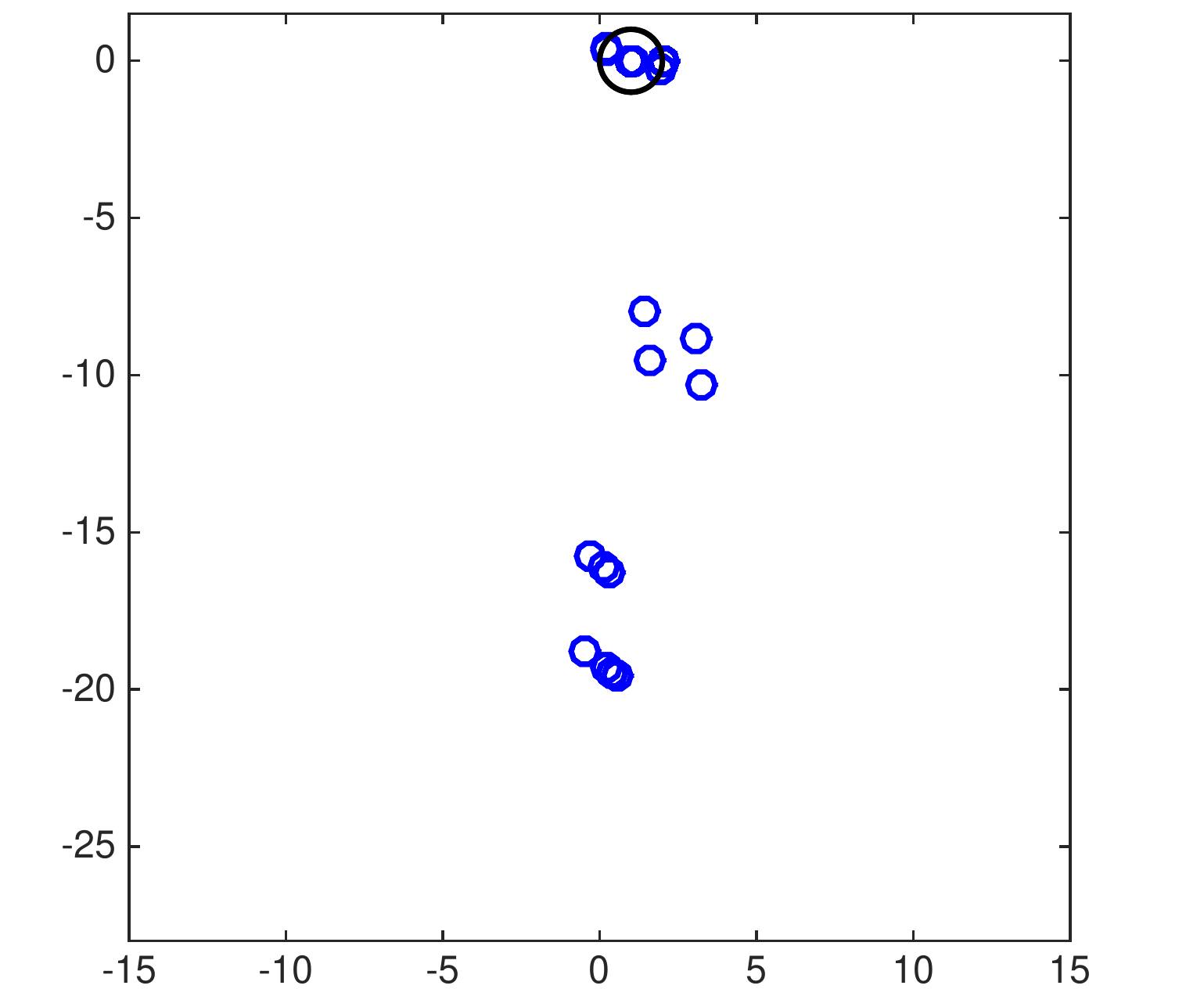}} \quad
\subfloat[][$k=3$]
  {\includegraphics[width=.45\textwidth]{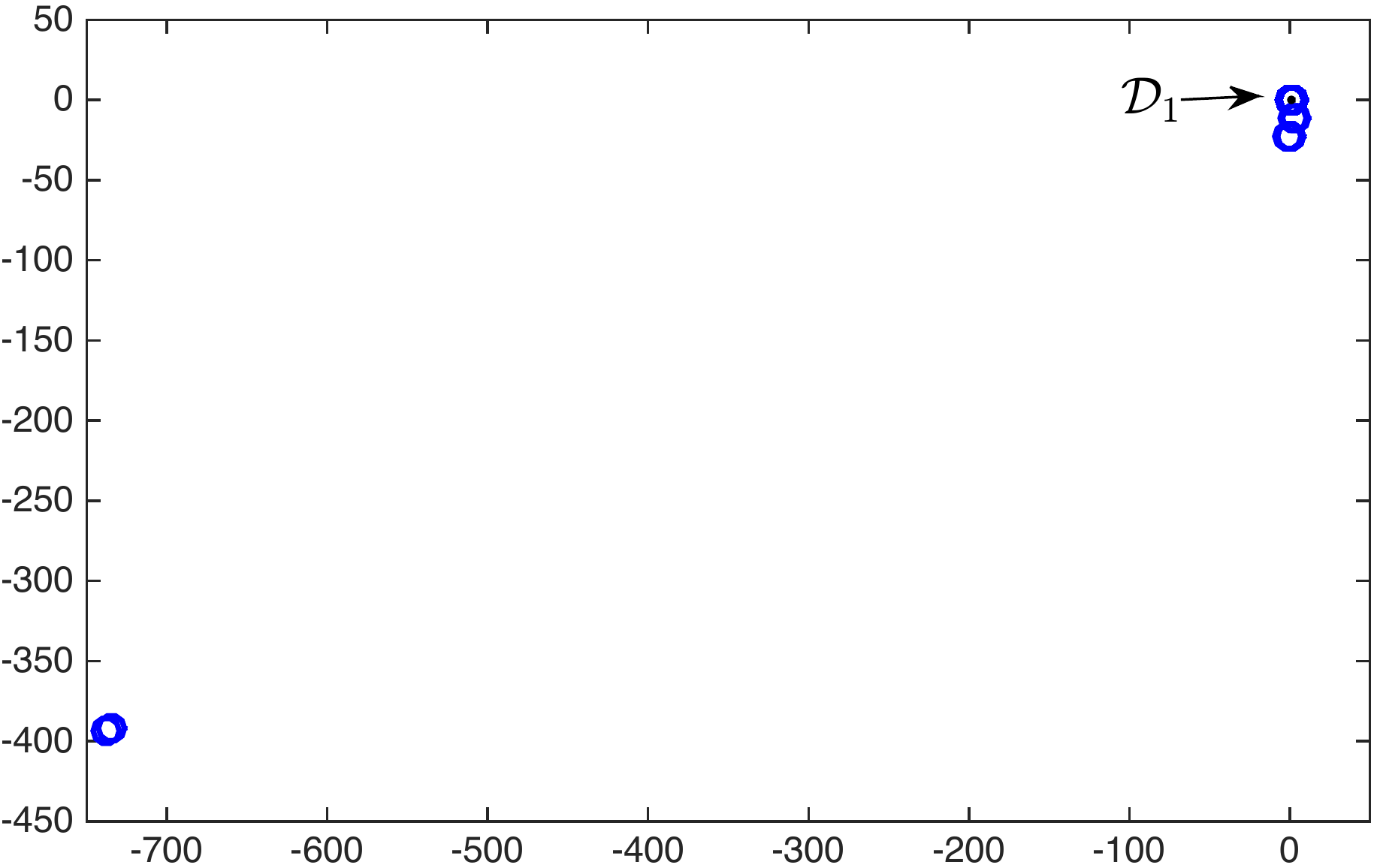}} 
\caption{Influence of the polynomial degree $r=k+1$ on the spectrum of the OAS-preconditioned matrix for $\omega=\omega_2$, $N_{\text{sub}}=2$, $\delta_\text{ovr}=2h$.}
\label{fig:spectraASvariak}
\end{figure}

\begin{figure}
\centering
\subfloat[][$\omega=\omega_1$]
  {\includegraphics[width=.45\textwidth]{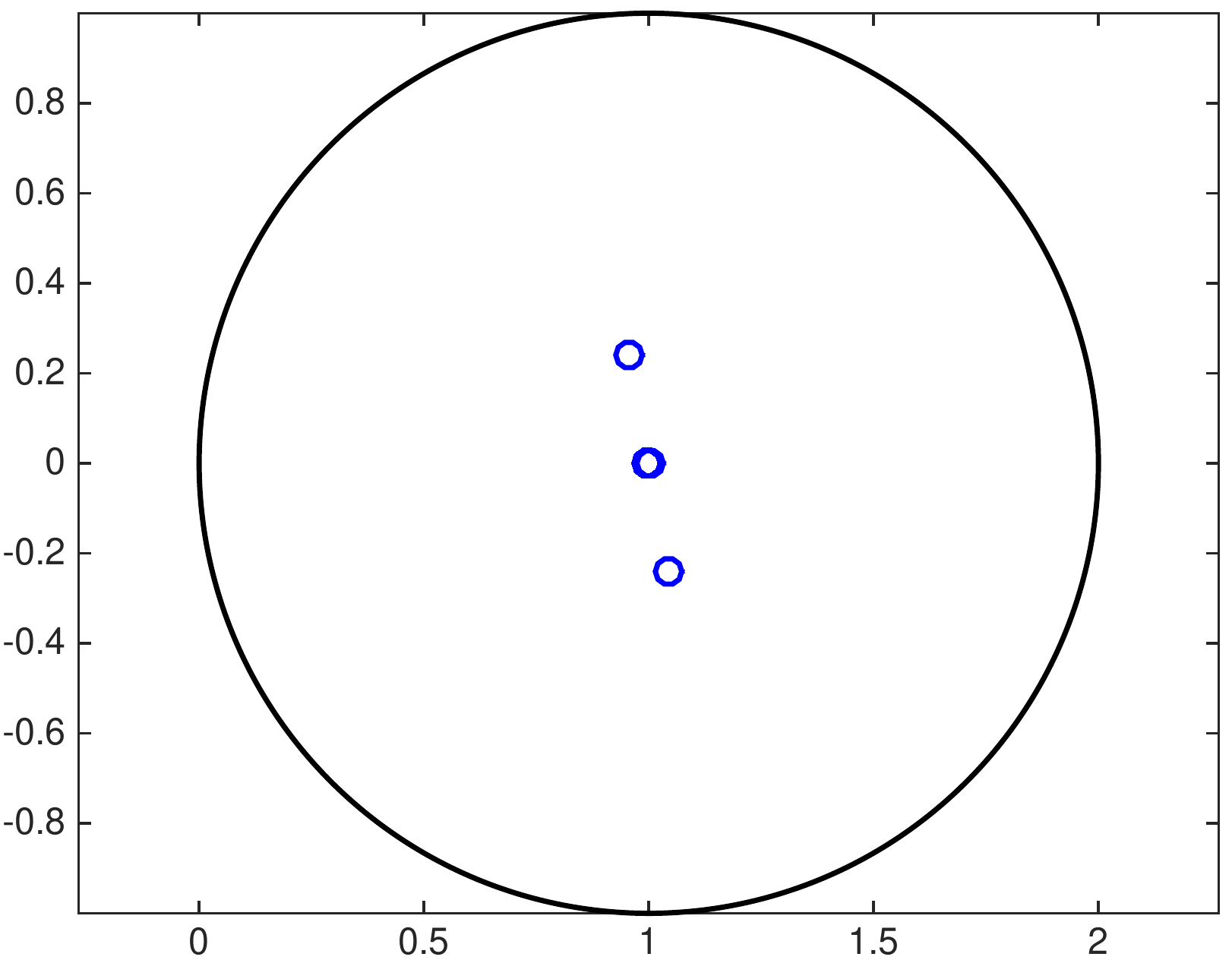}} \quad
\subfloat[][$\omega=\omega_3$]
  {\includegraphics[width=.45\textwidth]{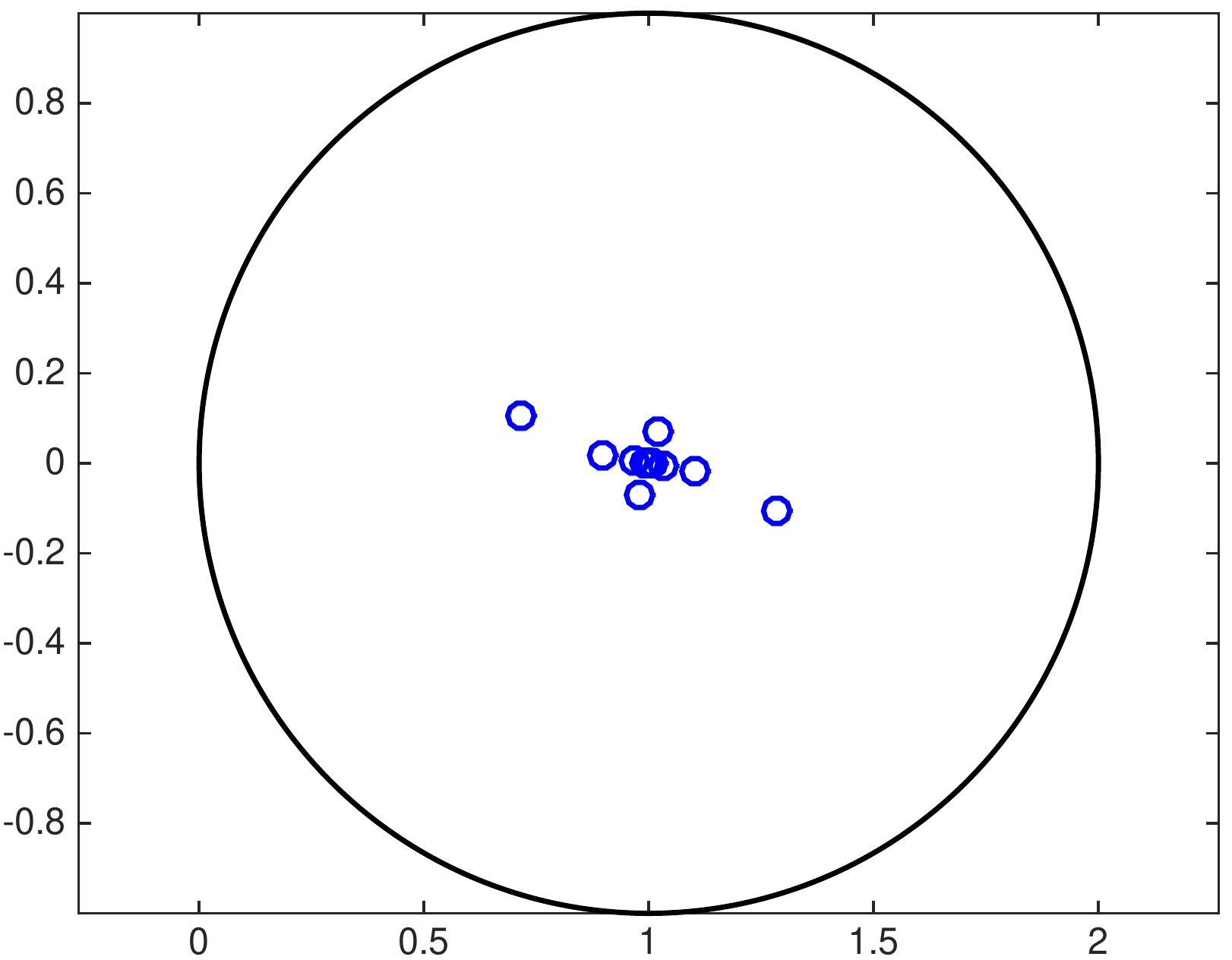}}
\caption{Influence of the angular frequency $\omega$ on the spectrum of the ORAS-preconditioned matrix for $k=2$, $N_{\text{sub}}=2$,
 $\delta_\text{ovr}=2h$.}
\label{fig:spectraRASvariaw}
\end{figure}

\begin{figure}
\centering
\subfloat[][$\omega=\omega_1$]
  {\includegraphics[width=.45\textwidth]{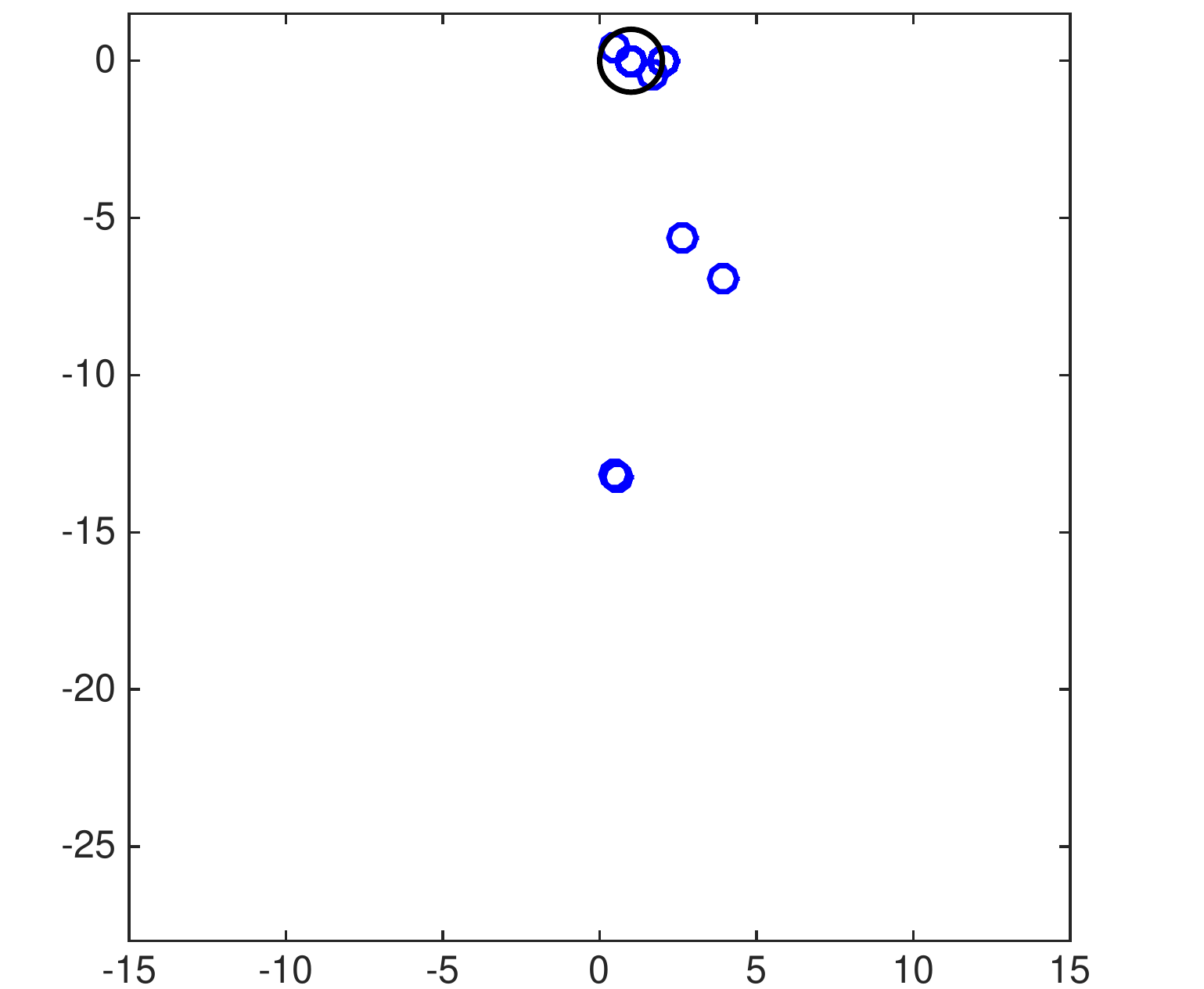}} \quad
\subfloat[][$\omega=\omega_3$]
  {\includegraphics[width=.45\textwidth]{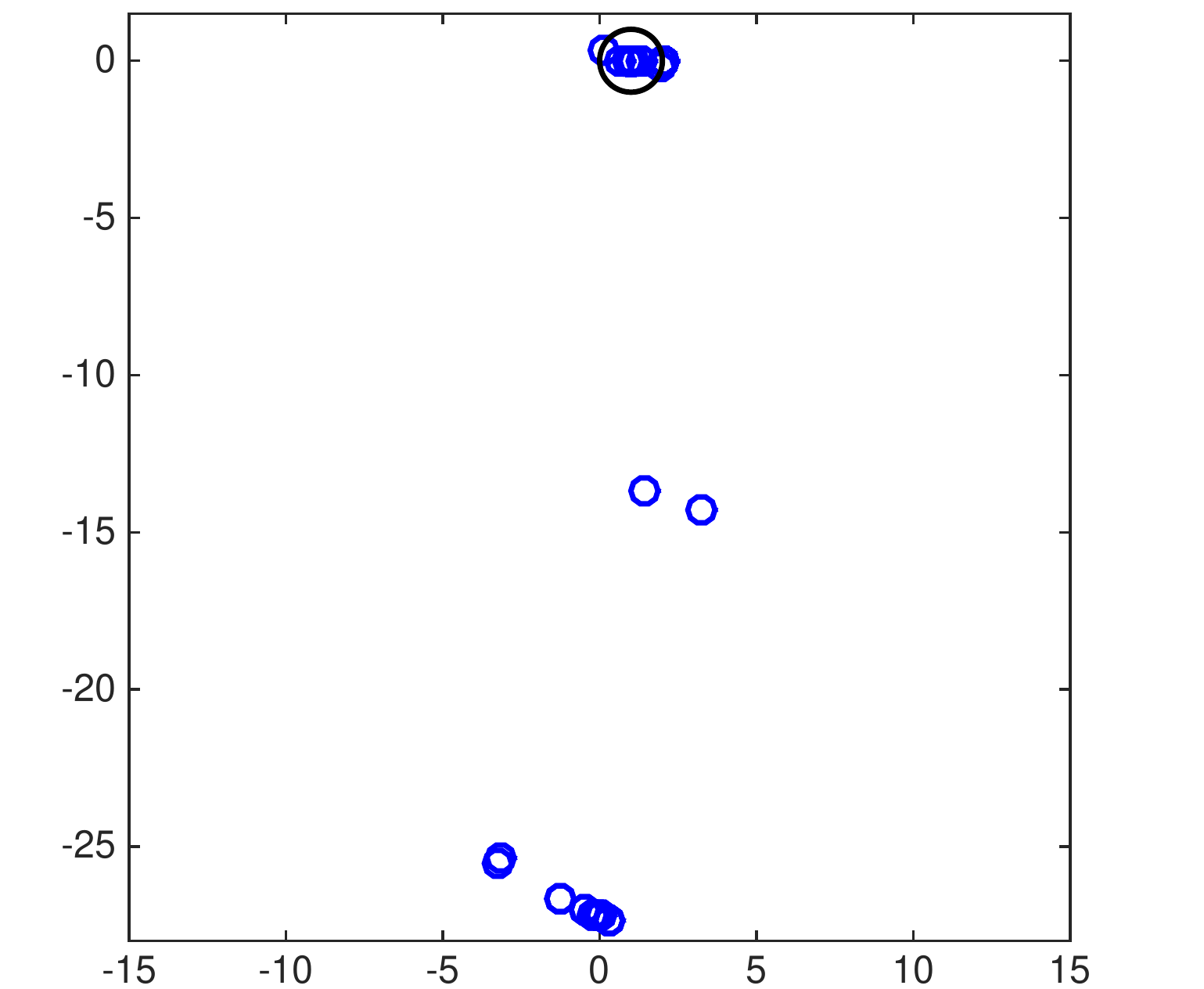}}
\caption{Influence of the angular frequency $\omega$ on the spectrum of the OAS-preconditioned matrix for $k=2$, $N_{\text{sub}}=2$,
 $\delta_\text{ovr}=2h$.}
\label{fig:spectraASvariaw}
\end{figure}

\begin{figure}
\centering
\subfloat[][$N_{\text{sub}}=4$]
  {\includegraphics[width=.45\textwidth]{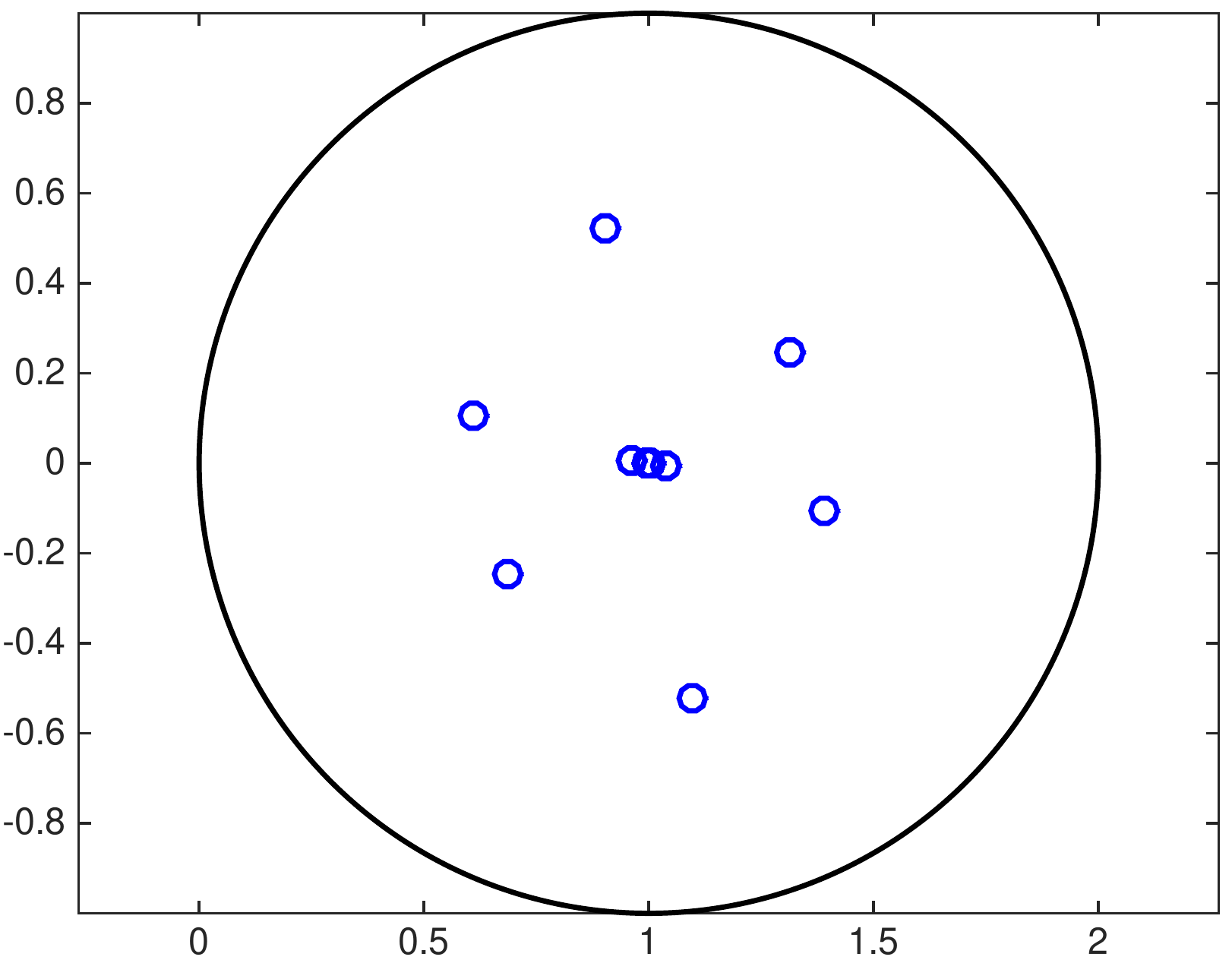}} \quad
\subfloat[][$N_{\text{sub}}=8$]
  {\includegraphics[width=.45\textwidth]{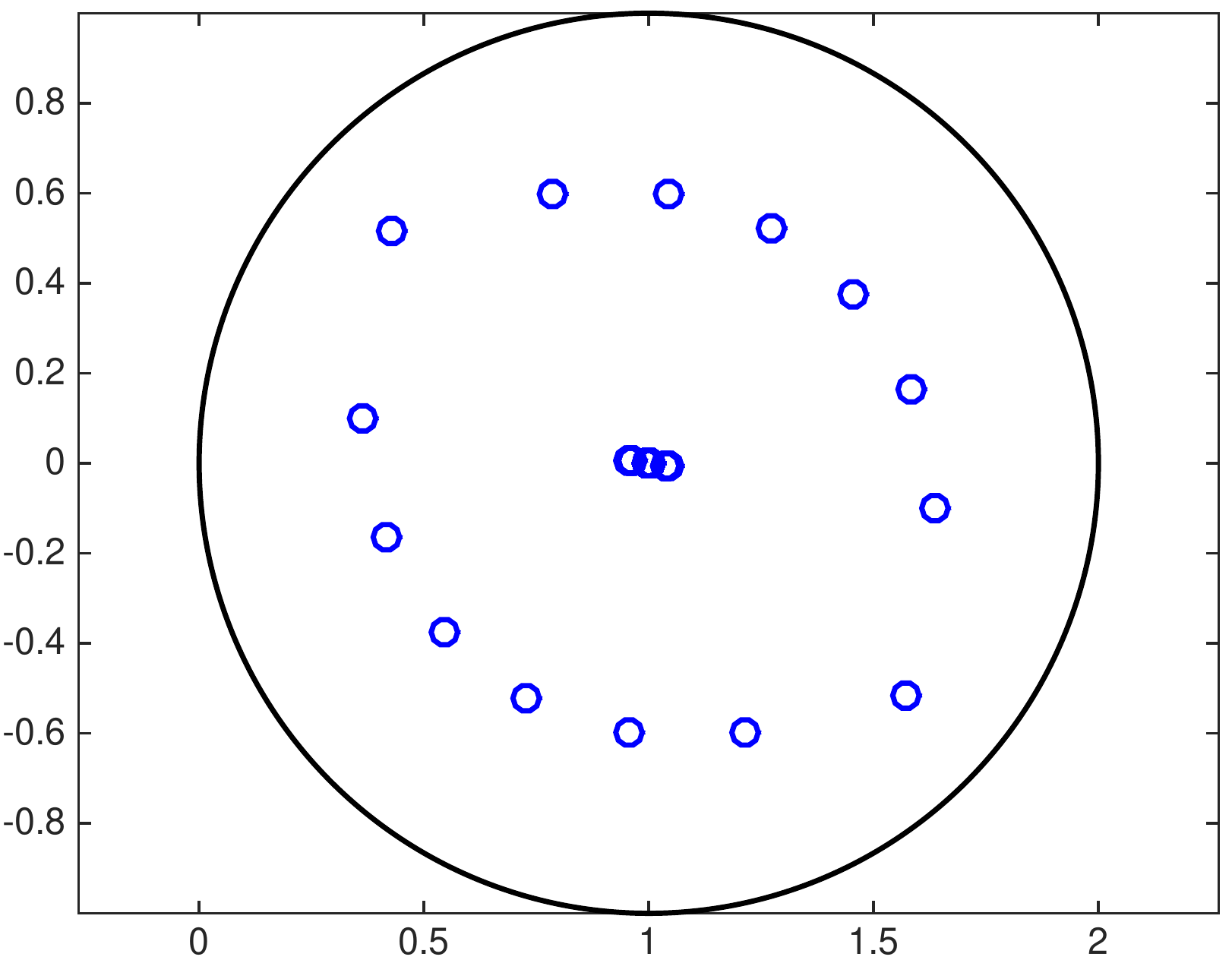}}
\caption{Influence of the number of subdomains $N_{\text{sub}}$ on the spectrum of the ORAS-preconditioned matrix for $k=2$, $\omega=\omega_2$, $\delta_\text{ovr}=2h$.}
\label{fig:spectraRASvariaNsub}
\end{figure}

\begin{figure}
\centering
\subfloat[][$N_{\text{sub}}=4$]
  {\includegraphics[width=.45\textwidth]{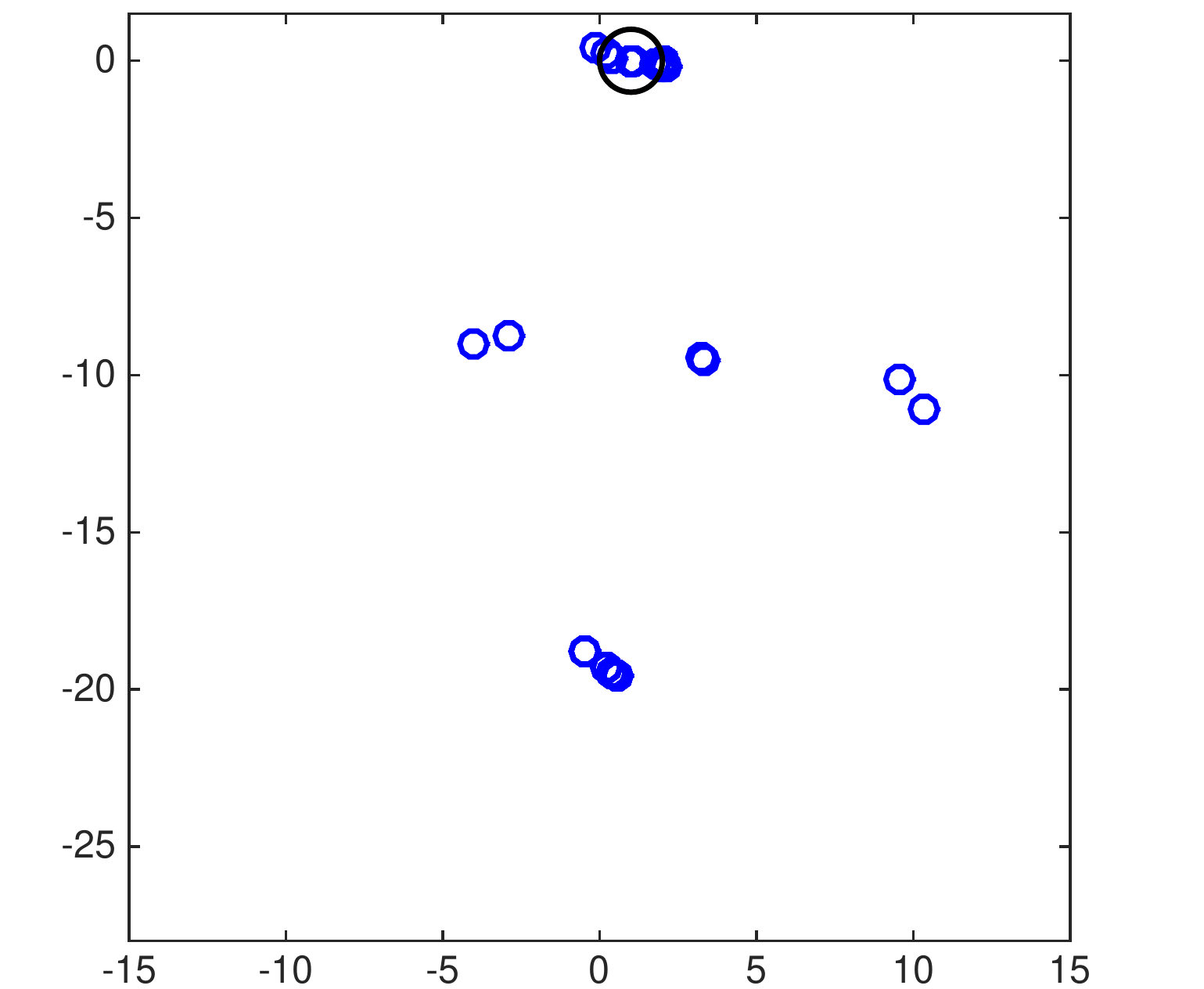}} \quad
\subfloat[][$N_{\text{sub}}=8$]
  {\includegraphics[width=.45\textwidth]{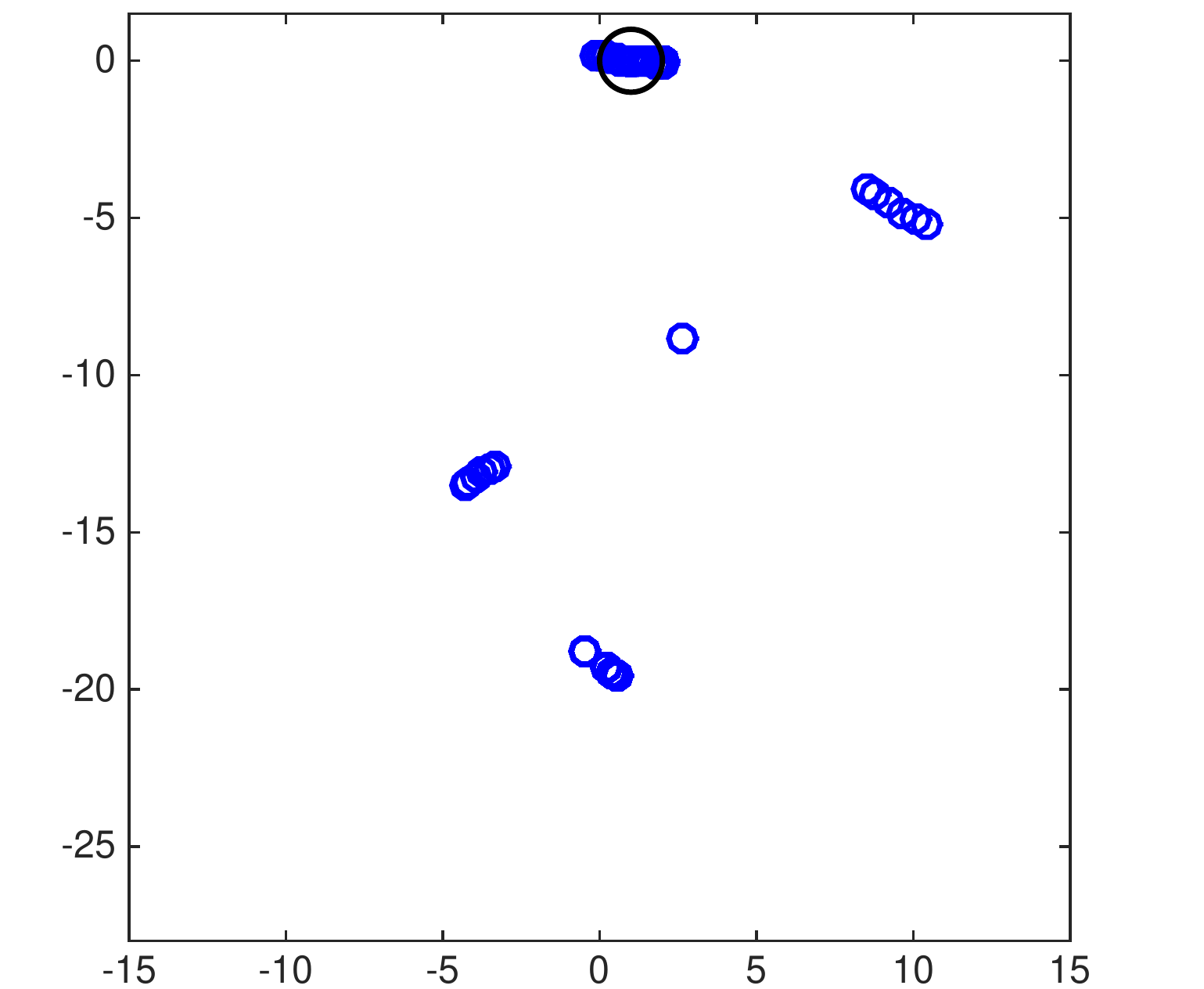}}
\caption{Influence of the number of subdomains $N_{\text{sub}}$ on the spectrum of the OAS-preconditioned matrix for $k=2$, $\omega=\omega_2$, $\delta_\text{ovr}=2h$.}
\label{fig:spectraASvariaNsub}
\end{figure}

\begin{figure}
\centering
\subfloat[][$\delta_\text{ovr}=2h$]
  {\includegraphics[width=.45\textwidth]{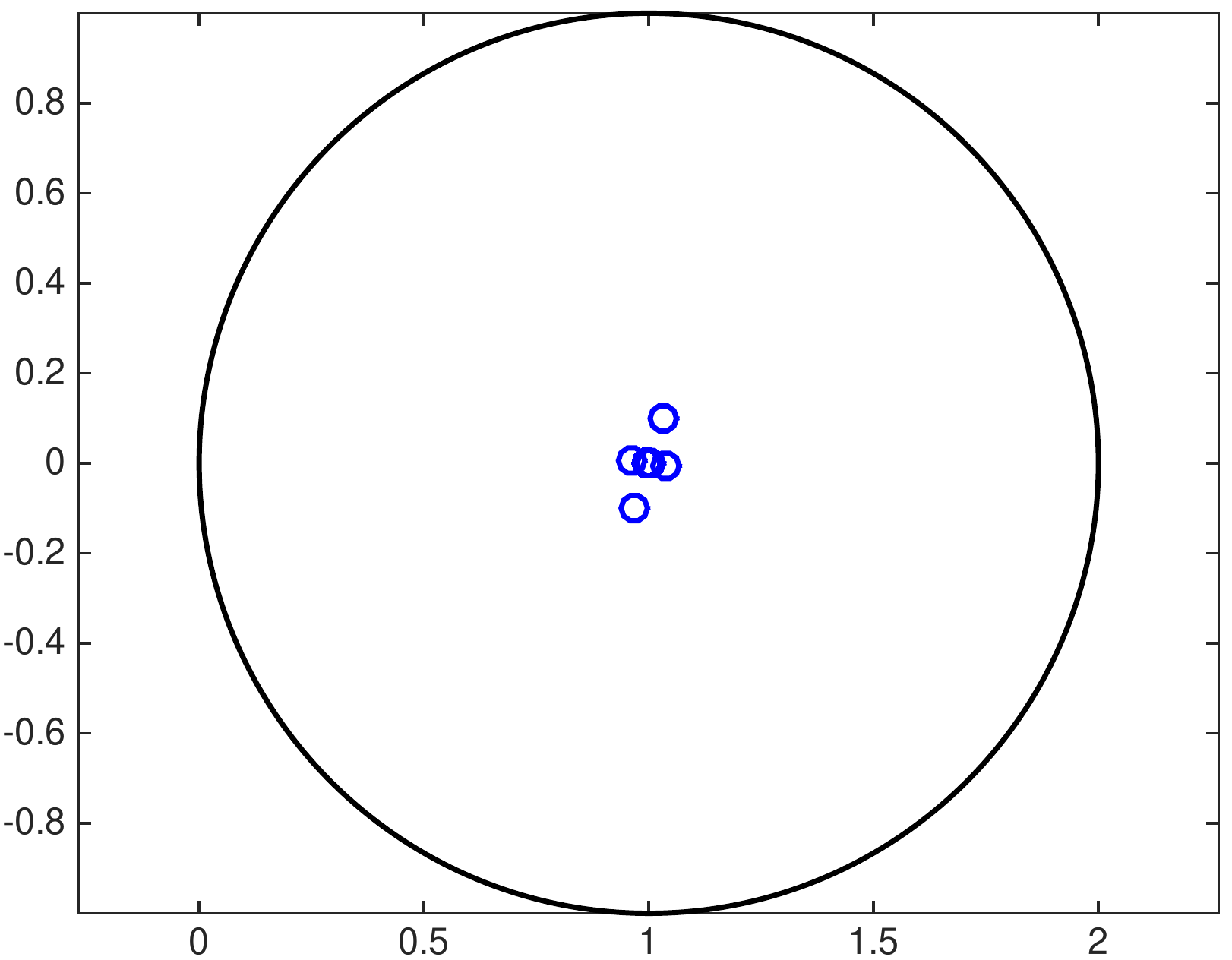}} \quad
\subfloat[][$\delta_\text{ovr}=4h$]
  {\includegraphics[width=.45\textwidth]{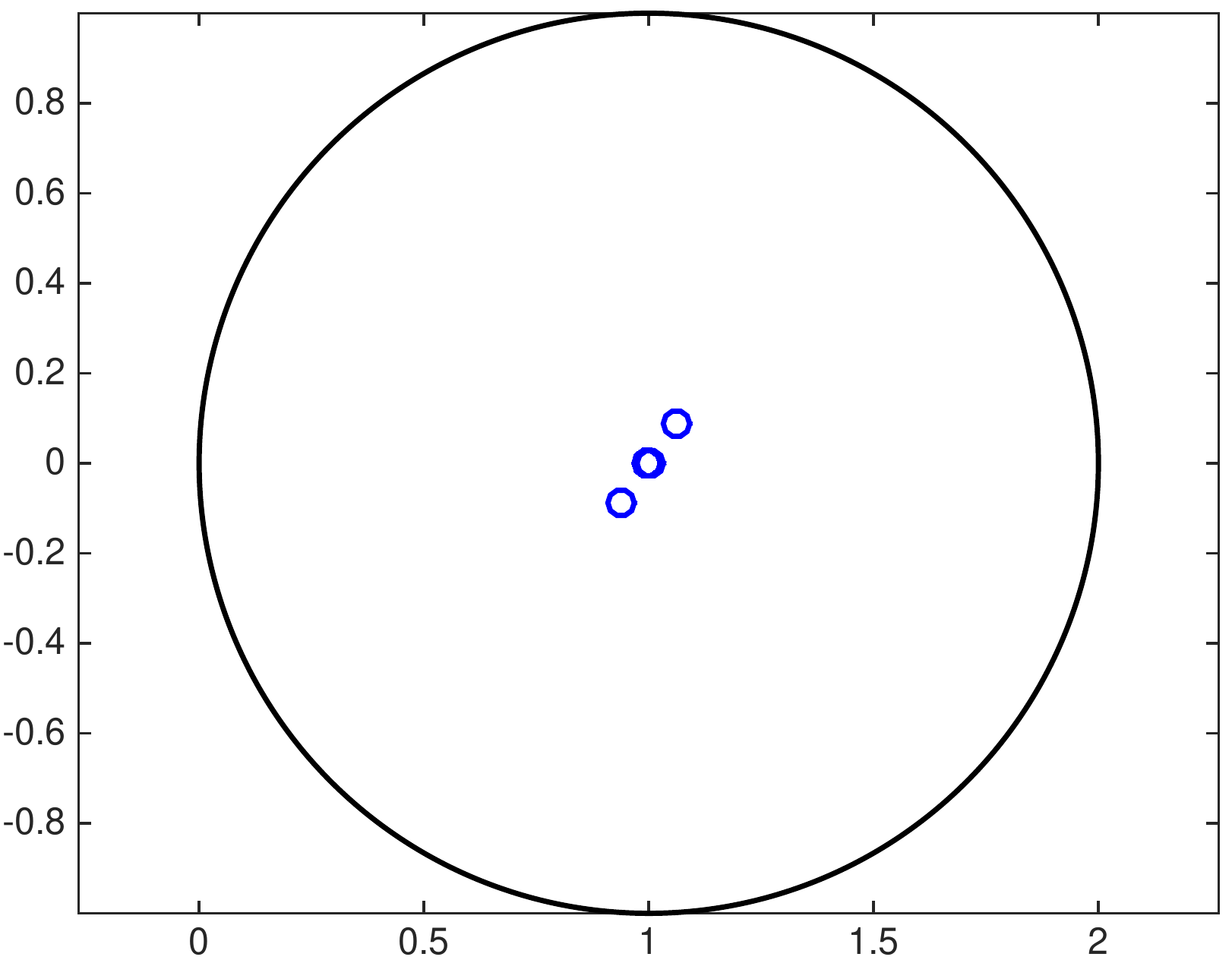}}
\caption{Influence of the overlap size $\delta_\text{ovr}$ on the spectrum of the ORAS-preconditioned matrix for $k=2$, $\omega=\omega_2$,
 $N_{\text{sub}}=2$.}
\label{fig:spectraRASvariaovr}
\end{figure}

\begin{figure}
\centering
\subfloat[][$\delta_\text{ovr}=2h$]
  {\includegraphics[width=.45\textwidth]{ASspectrumk2nSub2ovr1w2hb}} \quad
\subfloat[][$\delta_\text{ovr}=4h$]
  {\includegraphics[width=.45\textwidth]{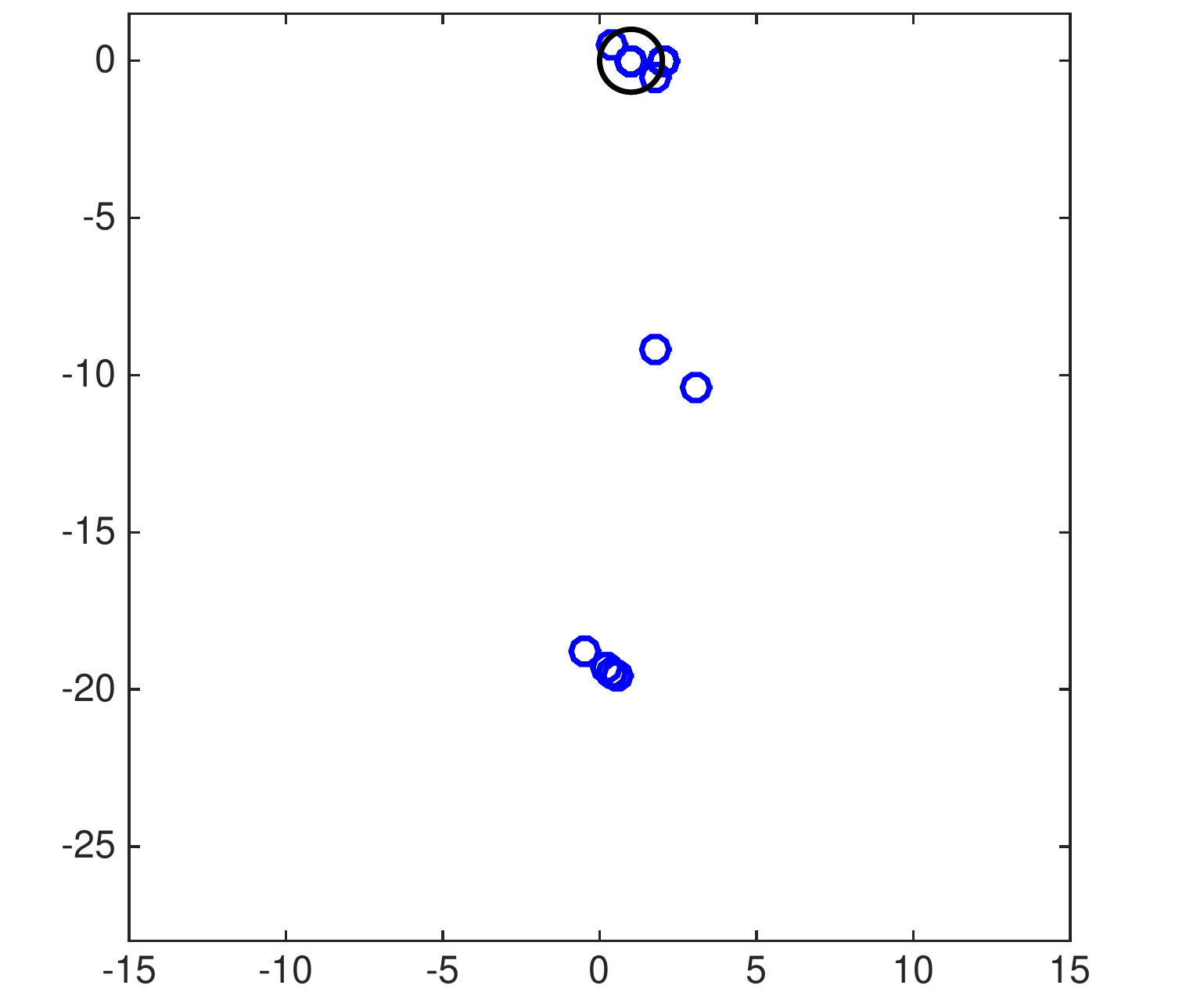}}
\caption{Influence of the overlap size $\delta_\text{ovr}$ on the spectrum of the OAS-preconditioned matrix for $k=2$, $\omega=\omega_2$,
 $N_{\text{sub}}=2$.}
\label{fig:spectraASvariaovr}
\end{figure}

To study the convergence of GMRES preconditioned by ORAS or OAS we vary
first the polynomial degree $r=k+1$ (Table~\ref{tab:variak}, Figs.~\ref{fig:spectraRASvariak}--\ref{fig:spectraASvariak}), 
then the angular frequency $\omega$ (Table~\ref{tab:variaw}, Figs.~\ref{fig:spectraRASvariaw}--\ref{fig:spectraASvariaw}), 
the number of subdomains $N_{\text{sub}}$ (Table~\ref{tab:variaNsub}, Figs.~\ref{fig:spectraRASvariaNsub}--\ref{fig:spectraASvariaNsub}) 
and finally the overlap size $\delta_\text{ovr}$ (Table~\ref{tab:variaovr}, Figs.~\ref{fig:spectraRASvariaovr}--\ref{fig:spectraASvariaovr}). 
Here, $\delta_\text{ovr} = 1h, 2h, 4h$ means that we consider a total overlap between two subdomains of $1,2,4$ mesh triangles along the horizontal direction (see Fig.~\ref{fig:decomp}). 

In Tables~\ref{tab:variak}--\ref{tab:variaovr}, $N_{\text{dofs}}$ is the total number of degrees of freedom, 
$N_{\text{iterNp}}$ is the number of iterations necessary to attain the prescribed convergence for GMRES without any preconditioner, 
and $N_{\text{iter}}$ is the number of iterations for GMRES preconditioned by ORAS (OAS).
Moreover, denoting by 
\[
{\cal D}_1=\{z\in\mathbb{C}: |z-z_0|<1\}
\]
the unit disk centered at $z_0=(1,0)$ in the complex plane, we measure also the maximum distance
to $(1,0)$ of the eigenvalues $\lambda$ of the preconditioned matrix, the number
of eigenvalues that have distance greater than $1$, and the number of
eigenvalues that have distance equal to $1$ (up to a tolerance of
$10^{-10}$). This information is useful to characterize the
convergence. Indeed, if $A$ is the matrix of the system to solve and $M^{-1}$ is the domain
decomposition preconditioner, then $I-M^{-1}A$ is the iteration matrix
of the domain decomposition method used as an iterative solver. 
So, a measure of the convergence of the domain decomposition solver
would be to check whether the eigenvalues of the preconditioned matrix $M^{-1}A$ are contained in ${\cal D}_1$. 
When the domain decomposition method is used, like here, as a preconditioner, the distribution of the spectrum remains a good indicator of the convergence.
Note that the matrix of the linear
system doesn't change when $N_{\text{sub}}$ or $\delta_\text{ovr}$
vary, therefore in Tables~\ref{tab:variaNsub}--\ref{tab:variaovr} (where $k=2$, $\omega=\omega_2$) we
don't report $N_{\text{dofs}}=1806$ and $N_{\text{iterNp}}=1138$ again.
In all Tables~\ref{tab:variak}--\ref{tab:variaovr}, we don't mention the condition number of the preconditioned matrix: 
indeed, no convergence rate estimates in terms of the condition number of the matrix, 
as those we are used to with
the conjugate gradient method, are available for the GMRES method. 

Figs.~\ref{fig:spectraRASvariak}, \ref{fig:spectraRASvariaw}, \ref{fig:spectraRASvariaNsub}, \ref{fig:spectraRASvariaovr}, respectively Figs.~\ref{fig:spectraASvariak}, \ref{fig:spectraASvariaw}, \ref{fig:spectraASvariaNsub}, \ref{fig:spectraASvariaovr}, show the whole spectrum in the complex plane of the matrix preconditioned by ORAS, respectively by OAS (note that many eigenvalues are multiple), together with ${\partial{\cal D}_1}$. 





\begin{figure}
\centering
\includegraphics[width=0.45\textwidth]{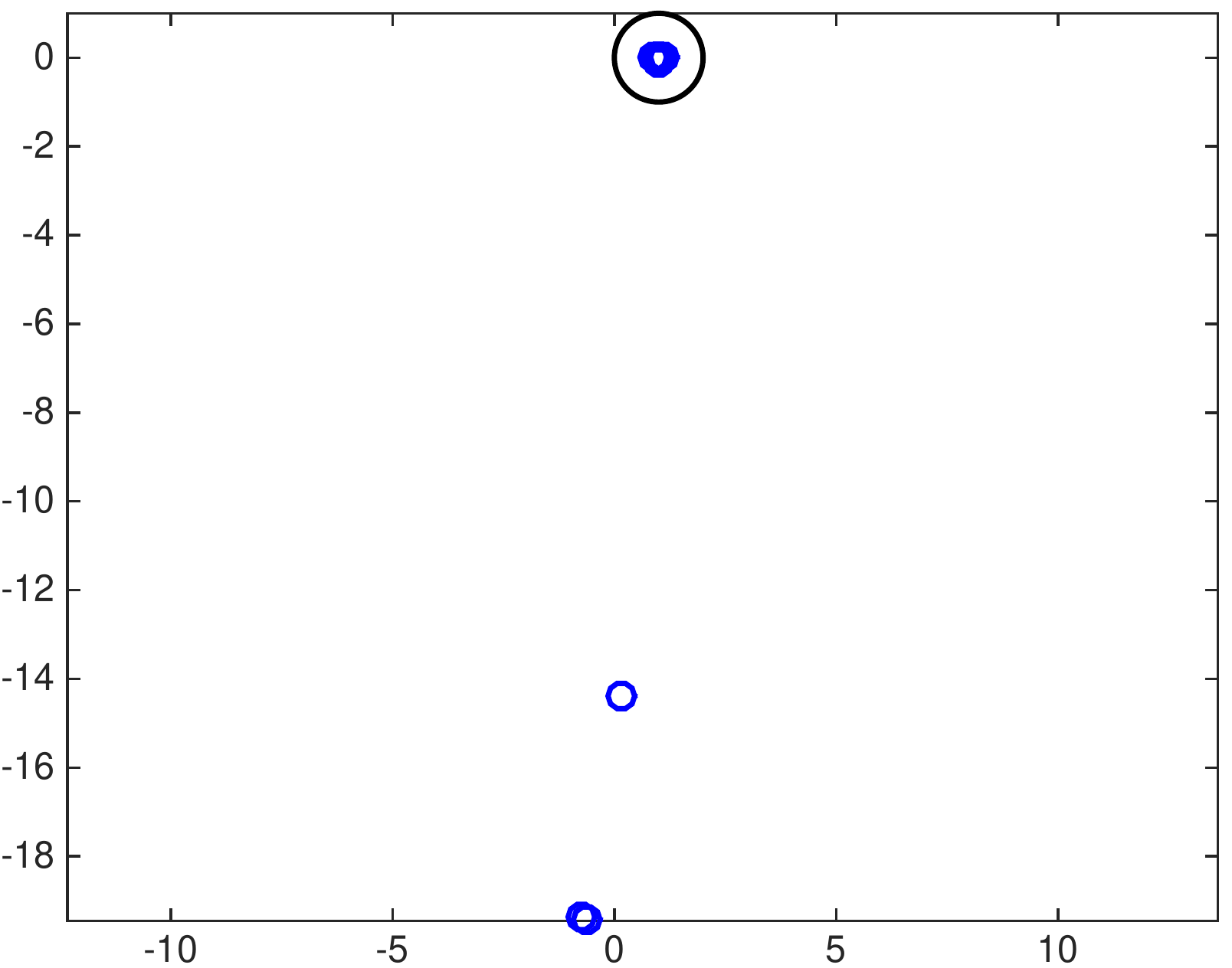}
\caption{The spectrum of the ORAS-preconditioned matrix for $k=2$, $\omega=\omega_2$,
 $N_{\text{sub}}=2$, $\delta_\text{ovr}=1h$.}
\label{fig:spectrumRASovr1h}
\end{figure}

Looking at the tables and figures, we can see that the non preconditioned GMRES method is very slow, and
the ORAS preconditioner gives much faster convergence than the OAS
preconditioner. As expected, the convergence becomes slower when $\omega$ or $N_{\text{sub}}$ increase, or when $\delta_\text{ovr}$ decreases. In these tests, when varying $k$ (which gives the polynomial degree $r=k+1$ of the FE basis functions), the number of iterations for convergence using the ORAS preconditioner is equal to $5$ for $k=0$ and then it stays equal to $6$ for $k>0$; this is reflected by the corresponding spectra in Fig.~\ref{fig:spectraRASvariak}, which indeed remain quite similar when $k$ varies.

Note also that, when using the ORAS preconditioner, for $2$ subdomains the spectrum is always well clustered inside the unit disk, except for the case with $\delta_\text{ovr}=1h$ (see Fig.~\ref{fig:spectrumRASovr1h}), in which $3$ eigenvalues are outside with distances from $(1,0)$ equal to $19.5, 19.4, 14.4$. This case $\delta_\text{ovr}=1h$ corresponds to adding a layer of simplices just to one of the two non overlapping subdomains to obtain the overlapping decomposition; hence it appears necessary to add at least one layer from both subdomains.
Then we see that for $4$ and $8$ subdomains 
the spectrum becomes less well clustered.
With the OAS preconditioner there are always eigenvalues outside the unit disk. 
For all the considered cases, we see that the less clustered the spectrum, the slower the convergence.

\subsection{Results for the three-dimensional problem}

We complete the presentation showing some results for the full 3d simulation, for a waveguide of dimensions $\mathtt{c} = 0.1004$\,m, $\mathtt{b} = 0.00508$\,m, and $\mathtt{a} = 0.01016$\,m.
The physical parameters are: $\varepsilon = 8.85\cdot10^{-12}$\,F\,m$^{-1}$, $\mu = 1.26\cdot10^{-6}$\,H\,m$^{-1}$ and $\sigma = 0.15$\,S\,m$^{-1}$ or $\sigma = 0$\,S\,m$^{-1}$.
We take a stripwise subdomains decomposition along the wave propagation, with $\delta_\text{ovr}=2h$; however, note that in FreeFem++ very general subdomains decompositions can be considered.

In 3d, if $\sigma=0$ there is an exact solution given by the Transverse Electric (TE) modes:
\begin{align*}
& E_x ^{TE} = 0, \\
& E_y ^{TE} = -\mathcal{C} \frac{m\pi}{\mathtt{a}} \sin\left(\frac{m\pi z}{\mathtt{a}}\right)\cos\left(\frac{n\pi y}{\mathtt{b}}\right) e^{-\mathtt{i}\beta x}, \\
& E_z ^{TE} = \mathcal{C} \frac{n\pi}{\mathtt{b}}\cos\left(\frac{m\pi z}{\mathtt{a}}\right)\sin\left(\frac{n\pi y}{\mathtt{b}}\right) e^{-\mathtt{i}\beta x}, 
\qquad m,n \in \mathbb{N}.
\end{align*}
The real constant $\beta$ is linked to the waveguide dimensions $\mathtt{a},\mathtt{b}$ by the so called dispersion relation $\left(\frac{m\pi}{\mathtt{a}}\right)^2+\left(\frac{n\pi}{\mathtt{b}}\right)^2 = \tilde \omega^2-\beta^2$, and we choose $\mathcal{C} = \mathtt{i} \omega \mu / (\tilde \omega^2-\beta^2)$. 
The field $\mathbf{E}^{TE}$ satisfies the metallic boundary conditions on $\Gamma_\text{w}$ and the impedance boundary conditions on $\Gamma_{\text{in}}$, $\Gamma_{\text{out}}$ with parameter $\eta=\beta$ and $\mathbf{g}^{\text{in}} = (\mathtt{i}\beta + \mathtt{i} \beta)\mathbf{E}^{TE} = 2\mathtt{i}\beta \mathbf{E}^{TE}$ and $\mathbf{g}^{\text{out}} = (-\mathtt{i}\beta + \mathtt{i} \beta)\mathbf{E}^{TE} = \mathbf{0}$. 

Since the propagation constant in 3d is $\beta$ and no more $\tilde\omega$, we compute the mesh size $h$ using the relation $h^2\cdot\beta^3=1$, taking $\beta = \omega_\beta \sqrt{\mu \varepsilon}$, with $\omega_\beta = 32$\,GHz.
Then the dispersion relation gives $\tilde\omega = \sqrt{ \beta^2 + ({m\pi}/{\mathtt{a}})^2+({n\pi}/{\mathtt{b}} )^2 }$ (where we choose $m=1,n=0$), and we get $\omega = \tilde\omega/ \sqrt{\mu \varepsilon}$.

Again, the linear system is solved with preconditioned GMRES, with a stopping criterion based on the relative residual and a tolerance of $10^{-6}$, starting with a random initial guess.
To apply the preconditioner, the local problems in each subdomain of matrices $A_s$ are solved with the direct solver MUMPS \cite{amestoy:2001:fully}.

\begin{table}
\centering
$ \begin{array}{crrccrr}
\toprule
k & N_{\text{dofs}} & N_{\text{iter}} & & N_{\text{sub}} & N_{\text{dofs}} & N_{\text{iter}}  \\
\midrule
0	& 62283   & 8(40)   & & 2 & 324654 & 8(70) \\
1	& 324654 & 8(70)   & & 4 & 324654 & 11(106) \\
2	& 930969 & 8(99) & & 8 & 324654 & 17(168)\\
\bottomrule
\end{array}$
\caption{Results in 3d, $\sigma = 0.15$\,S\,m$^{-1}$: influence of the polynomial degree $r=k+1$ (for $N_{\text{sub}}=2$), and of the number of subdomains $N_{\text{sub}}$ (for $k=1$), on the convergence of ORAS(OAS) preconditioner ($\beta = \omega_\beta \sqrt{\mu \varepsilon}$ with $\omega_\beta = 32$\,GHz, $\delta_\text{ovr}=2h$).}
\label{tab:3dresultsBiggersigma015}
\end{table}
\begin{table}
\centering
$ \begin{array}{crrccrr}
\toprule
k & N_{\text{dofs}} & N_{\text{iter}} & & N_{\text{sub}} & N_{\text{dofs}} & N_{\text{iter}}  \\
\midrule
0	& 62283   & 7(40)   & & 2 & 324654 & 8(67) \\
1	& 324654 & 8(67)   & & 4 & 324654 & 13(114) \\
2	& 930969 & 8(97) & & 8 & 324654 & 23(201)\\
\bottomrule
\end{array}$
\caption{Results in 3d, $\sigma = 0$\,S\,m$^{-1}$: influence of the polynomial degree $r=k+1$ (for $N_{\text{sub}}=2$), and of the number of subdomains $N_{\text{sub}}$ (for $k=1$), on the convergence of ORAS(OAS) preconditioner ($\beta = \omega_\beta \sqrt{\mu \varepsilon}$ with $\omega_\beta = 32$\,GHz, $\delta_\text{ovr}=2h$).}
\label{tab:3dresultsBiggersigma0}
\end{table}
In Tables~\ref{tab:3dresultsBiggersigma015}, \ref{tab:3dresultsBiggersigma0} we show the number of iterations for convergence, for the problem with $\sigma = 0.15$\,S\,m$^{-1}$ and $\sigma = 0$\,S\,m$^{-1}$ respectively, varying first the polynomial degree $r=k+1$ (for $N_{\text{sub}}=2$), and then the number of subdomains $N_{\text{sub}}$ (for $k=1$).
Like in the 2d case, the number of iterations using the ORAS preconditioner does not vary with the polynomial degree of the FE basis functions, while using the OAS preconditioner it varies and is much higher. Again, the convergence becomes slower when the number of subdomains increases, both with ORAS and OAS. We see that for more than $2$ subdomains the number of iterations for the non dissipative problem ($\sigma = 0$) is higher than for the problem with $\sigma = 0.15$\,S\,m$^{-1}$.

\begin{figure}
\centering
\subfloat
{\includegraphics[width=0.6\textwidth]{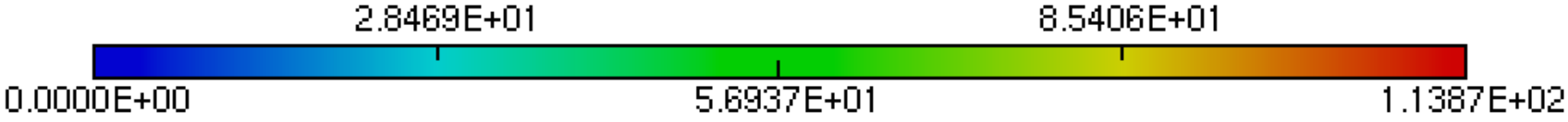}} \\
\subfloat
{\includegraphics[width=0.6\textwidth]{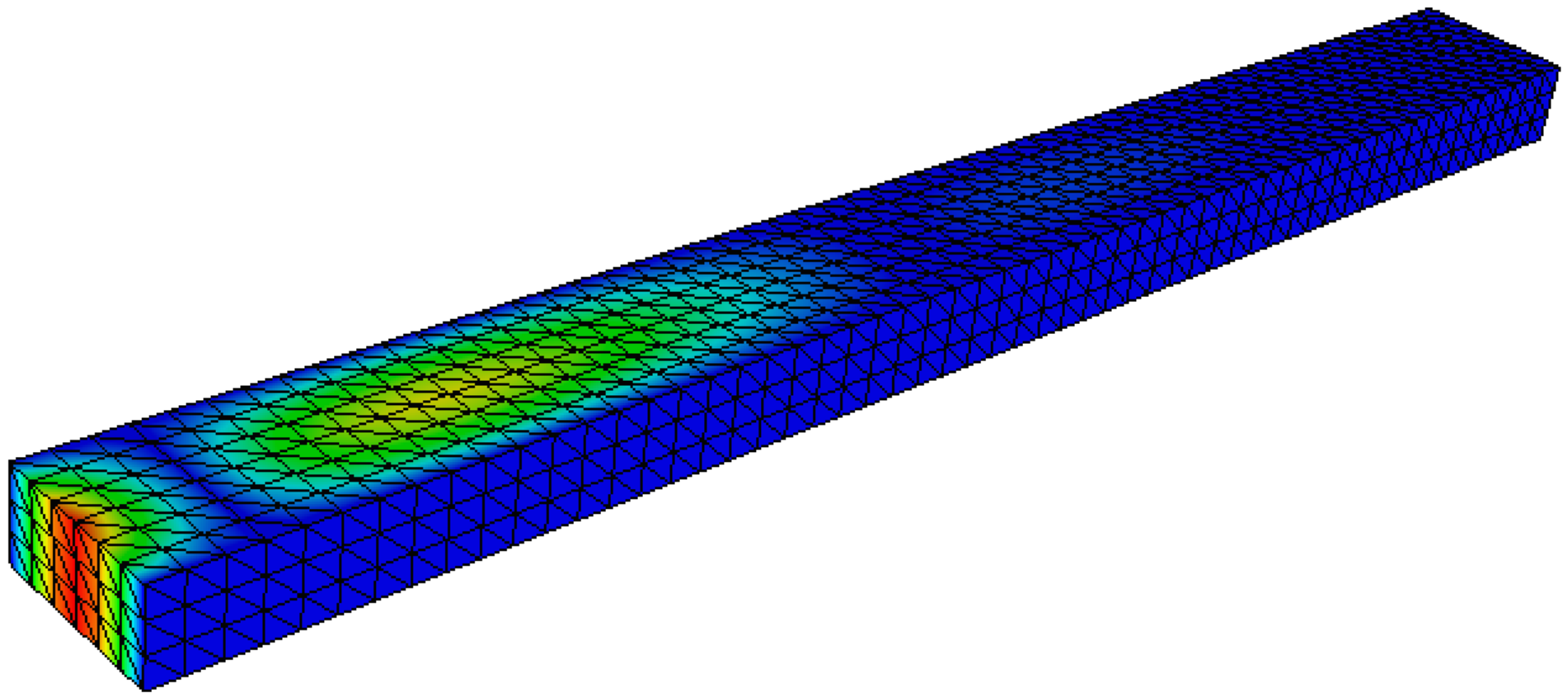}}\\
\subfloat
{\includegraphics[width=0.5\textwidth]{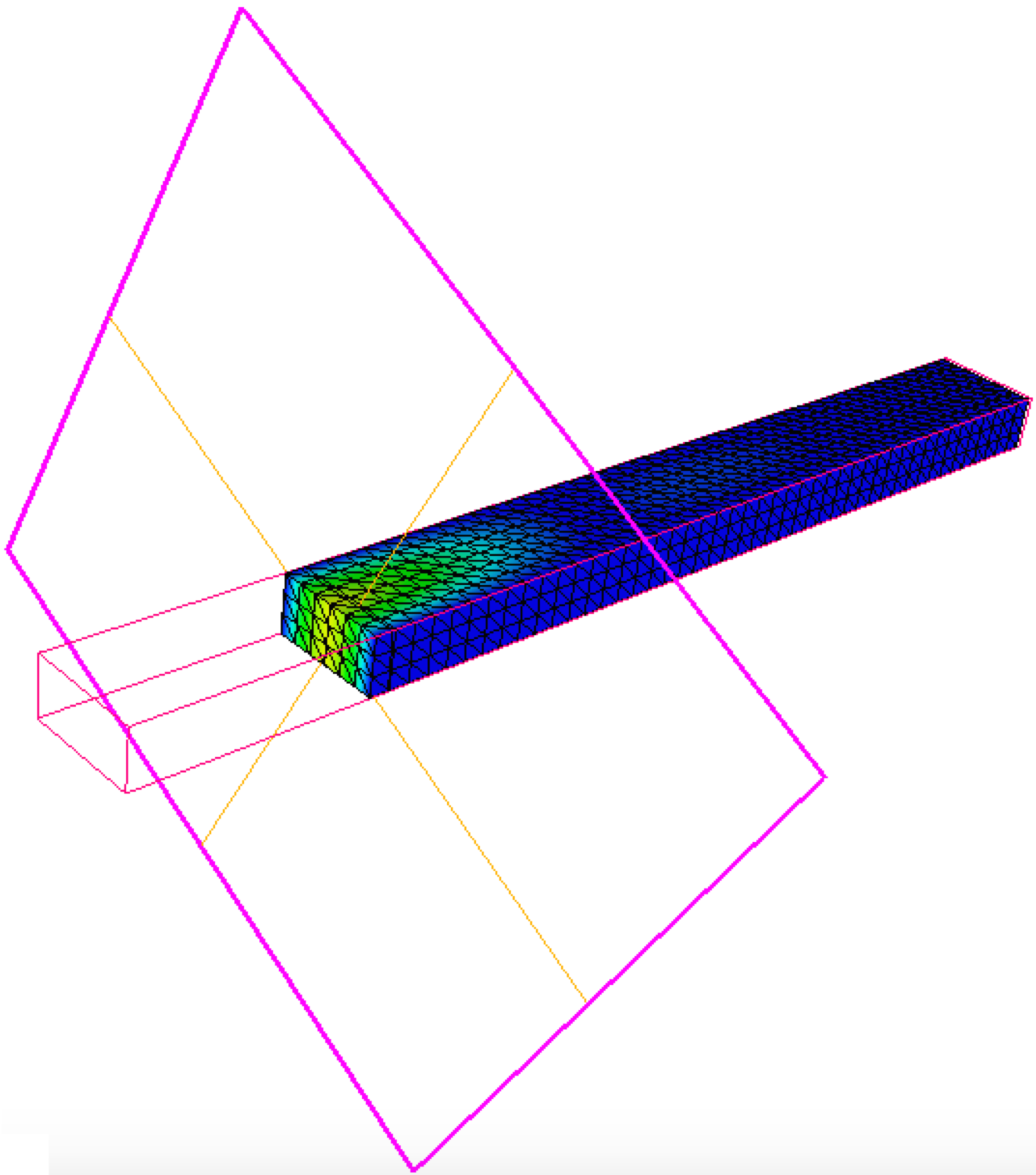}} \\
\subfloat
{\includegraphics[width=0.5\textwidth]{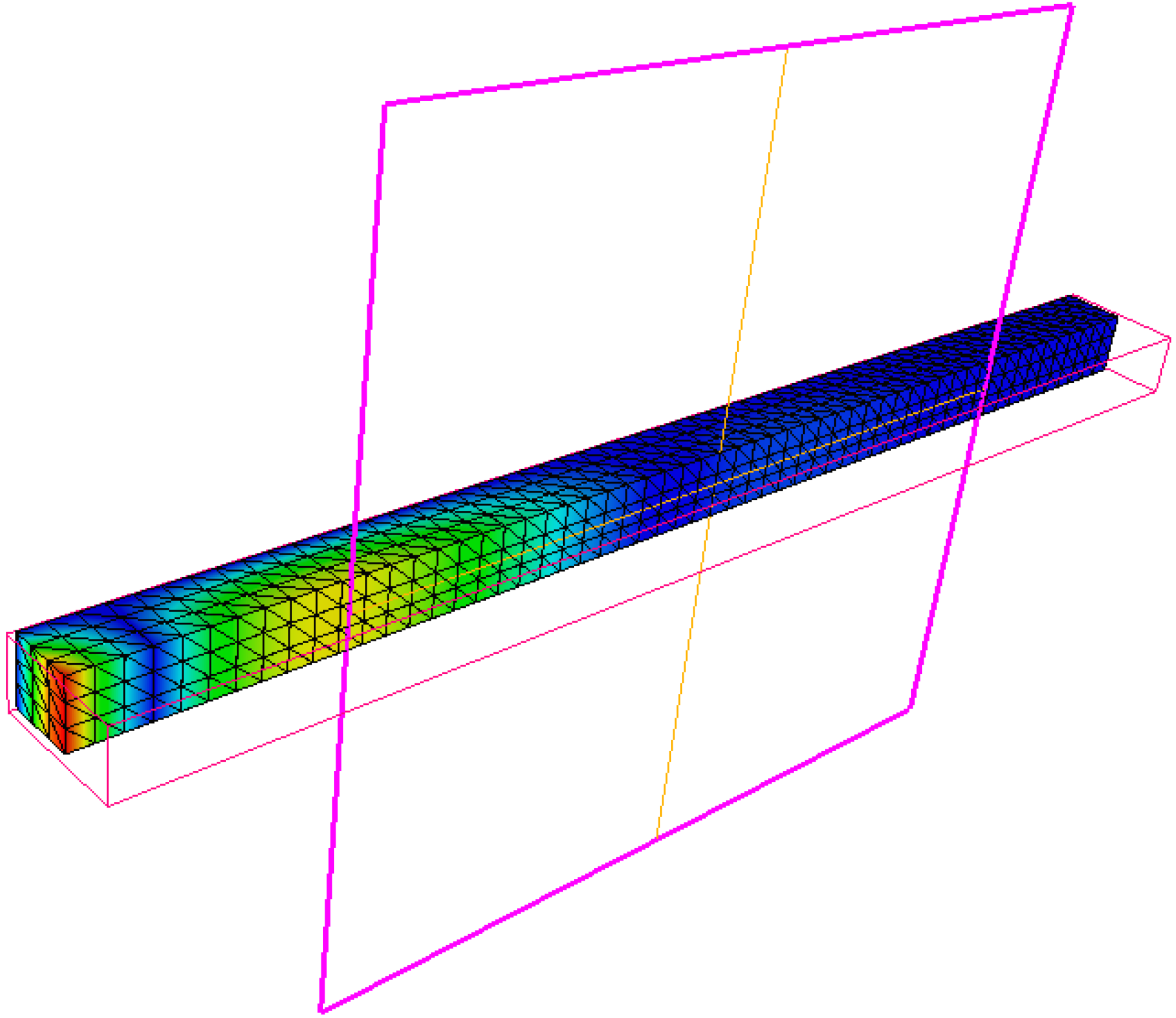}}
\caption{The norm of the real part of the solution for $\sigma = 0.15$\,S\,m$^{-1}$, with two sections of the waveguide.}
\label{fig:solution}
\end{figure}

In Figure~\ref{fig:solution} we plot the norm of the real part of the solution, which decreases as the wave propagates since there $\sigma = 0.15$\,S\,m$^{-1}$ is different from zero. 


\section{Conclusion}

We have adopted a friendly definition of high order edge elements generators and degrees of freedom: both in 2d and 3d their expression is rather simple, and the generators are strictly connected with the degrees of freedom.
Their presentation is enriched with illustrative examples, and an operational strategy of implementation of these elements is described in detail.
The elements in 3d of polynomial degree $1,2,3$ are available in FreeFem++ (since version 3.44), loading the plugin \texttt{load "Element\_Mixte3d"} and building the finite element space \texttt{fespace} with the keywords \texttt{Edge03d}, \texttt{Edge13d}, \texttt{Edge23d} respectively.

Numerical experiments have shown that Schwarz preconditioning significantly improves GMRES convergence for different values of physical and numerical parameters, and that the ORAS preconditioner always performs much better than the OAS preconditioner. 
Indeed, the only advantage of the OAS method is to preserve symmetry for symmetric problems: that is why it should be used only for symmetric positive definite matrices as a preconditioner for the conjugate gradient method.
Moreover, in all the considered test cases, the number of iterations for convergence using the ORAS preconditioner does not vary when the polynomial degree of the adopted high order finite elements increases.
We have also seen that it is necessary to take an overlap of at least one layer of simplices from \emph{both} subdomains of a neighbors pair. All these convergence qualities are reflected by the spectrum of the preconditioned matrix.

For higher order discretizations the computational cost per iteration grows since matrices become very large, therefore a parallel implementation as the one of HPDDM \cite{JolHecNat:2013:hpddm} (a high-performance unified framework for domain decomposition methods which is interfaced with FreeFem++) should be considered for large scale problems, as in \cite{BDRT:2017:EMF}. A two-level preconditioner via a coarse space correction should be designed for Maxwell's equations in order to fix the dependence on the number of subdomains or on the frequency of the iteration count.  
Note that in literature very few convergence theory results are available for domain decomposition methods applied to indefinite wave propagation problems. For the Helmholtz equation with absorption, rigorous convergence estimates have been recently presented in \cite{GraSpe:2017:DDH} for the two-level Additive Schwarz preconditioner, and their extension to the Maxwell case is work in progress.

\medskip
\noindent \textbf{Acknowledgement} 
This work was financed by the French National Research Agency (ANR) in the framework of the project MEDIMAX, ANR-13-MONU-0012.

\bibliographystyle{unsrt}
\bibliography{bdr}

\end{document}